\documentclass[11pt, draft]{amsart}

\usepackage{amsmath,amssymb,amsfonts,
}

\usepackage{diagrams}

\diagramstyle[centredisplay,dpi=600,nohug]
\newif\ifhavegoodboldmath\havegoodboldmathtrue
\newarrow{Tto}----{->}
\newarrow{Tinc}{\ifhavegoodboldmath bold\fi hook}---{->}
\newarrow{Tincdash}{\ifhavegoodboldmath bold\fi hook}{dash}{}{dash}{->}
\newarrow{Tonto}----{->>}
\newarrow{Tdash}{}{dash}{}{dash}{->}

\theoremstyle{plain}
\newtheorem{thm}{Theorem}[section]
\newtheorem{lem}[thm]{Lemma}
\newtheorem{cor}[thm]{Corollary}
\newtheorem{propo}[thm]{Proposition}

\theoremstyle{definition}
\newtheorem{defn}[thm]{Definition}
\newtheorem{ex}[thm]{Example}
\newtheorem{parraf}[thm]{}
\newtheorem*{ack}{Acknowledgments}

\theoremstyle{remark}
\newtheorem*{rem}{Remark}

\numberwithin{equation}{thm}

\newcommand{\HA}{\widehat{A}}
\newcommand{\HB}{\widehat{B}}
\newcommand{\hf}{\widehat{f}}
\newcommand{\hd}{\widehat{d}}

\newcommand{\BA}{\mathbb A}
\newcommand{\BD}{\mathbb D}
\newcommand{\NN}{\mathbb N}

\newcommand{\ZZ}{\mathbb Z}

\newcommand{\fm}{\mathfrak m}

\newcommand{\fp}{\mathfrak p}
\newcommand{\fq}{\mathfrak q}
\newcommand{\fr}{\mathfrak r}

\newcommand{\FS}{\mathfrak S}

\newcommand{\FU}{\mathfrak U}
\newcommand{\FV}{\mathfrak V}

\newcommand{\FX}{\mathfrak X}
\newcommand{\FY}{\mathfrak Y}
\newcommand{\FZ}{\mathfrak Z}

\newcommand{\sch}{\mathsf {Sch}}

\newcommand{\sfn}{\mathsf {NFS}}
\newcommand{\scha}{\mathsf {Sch}_{\mathsf {af}}}

\newcommand{\sfna}{\sfn_{\mathsf {af}}}

\newcommand{\CF}{\mathcal F}
\newcommand{\CG}{\mathcal G}
\newcommand{\CH}{\mathcal H}
\newcommand{\CI}{\mathcal I}
\newcommand{\CJ}{\mathcal J}
\newcommand{\CK}{\mathcal K}
\newcommand{\CL}{\mathcal L}

\newcommand{\CO}{\mathcal O}
\newcommand{\CP}{\mathcal P}
\newcommand{\CQ}{\mathcal Q}

\newcommand{\dirlim}[1]{\begin{array}[t]{c} \mathrm{lim}\\[-7.5 pt]
 {\longrightarrow} \\[-7.5 pt] {\scriptstyle {#1}} \end{array}}
\newcommand{\invlim}[1]{\begin{array}[t]{c} \mathrm{lim}\\[-7.5 pt]
 {\longleftarrow} \\[-7.5 pt] {\scriptstyle {#1}} \end{array}}

\newcommand{\lto}{\longrightarrow}
\newcommand{\xto}{\xrightarrow}

\newcommand{\epi}{\twoheadrightarrow}

\newcommand{\inc}{\hookrightarrow}

\newcommand{\imp}{\Rightarrow}
\newcommand{\dimp}{\Leftrightarrow}

\newcommand{\tr}{\triangle}
\newcommand{\tc}{\widehat{\otimes}}
\newcommand{\om}{\widehat{\Omega}}

\DeclareMathOperator{\Jac}{Jac}

\DeclareMathOperator{\rg}{rg}
\DeclareMathOperator{\htt}{ht}

\DeclareMathOperator{\spec}{Spec}

\DeclareMathOperator{\spf}{Spf}

\DeclareMathOperator{\Hom}{Hom}

\DeclareMathOperator{\ga}{\Gamma}

\DeclareMathOperator{\fD}{\mathfrak{D}}
\DeclareMathOperator{\dimtop}{dimtop}

\DeclareMathOperator{\coh}{\mathsf{Coh}}

\newcommand{\ie}{{\it i.e.} }

\begin{document}

\title{Local structure theorems for smooth maps of formal schemes}

\author[L. Alonso]{Leovigildo Alonso Tarr\'{\i}o}
\address{Departamento de \'Alxebra\\
Facultade de Matem\'a\-ticas\\
Universidade de Santiago de Compostela\\
E-15782  Santiago de Compostela, SPAIN}
\email{leoalonso@usc.es}
\urladdr{http://web.usc.es/\~{}lalonso/}

\author[A. Jerem\'{\i}as]{Ana Jerem\'{\i}as L\'opez}
\address{Departamento de \'Alxebra\\
Facultade de Matem\'a\-ticas\\
Universidade de Santiago de Compostela\\
E-15782  Santiago de Compostela, SPAIN}
\email{jeremias@usc.es}

\author[M. P\'erez]{Marta P\'erez Rodr\'{\i}guez}
\address{Departamento de Matem\'a\-ticas\\
Escola Superior de En\-xe\-\~ne\-r\'{\i}a Inform\'atica\\
Campus de Ourense, Univ. de Vigo\\
E-32004 Ou\-ren\-se, Spain}
\email{martapr@uvigo.es}

\thanks{This work was partially supported by Spain's MCyT and E.U.'s
FEDER research project MTM2005-05754.}


\begin{abstract}
We continue our study on infinitesimal lifting properties of maps between locally noetherian formal schemes started in \cite{AJP}. In this paper, we focus on some properties which arise specifically in the formal context. In this vein, we make a detailed study of the relationship between the infinitesimal lifting properties of a morphism of formal schemes and those of the corresponding maps of usual schemes associated to the directed systems that define the corresponding formal schemes. Among our main results, we obtain the characterization of completion morphisms as pseudo-closed immersions that are flat. Also, the local structure of smooth and \'etale  morphisms between locally noetherian formal schemes is described: the former factors locally as a completion morphism followed by a smooth adic morphism and the latter as a completion morphism followed by an \'etale adic morphism.
\end{abstract}



\maketitle

\hyphenation{pseu-do}

\maketitle

\tableofcontents

\section*{Introduction}
Formal schemes have always been present in the backstage of algebraic geometry but they were rarely studied in a systematic way after the foundational \cite[\S 10]{EGA1}. It has become more and more clear that the wide applicability of formal schemes in several areas of mathematics require such study. Let us cite a few of this applications. The construction of De Rham cohomology for a scheme $X$ of zero characteristic embeddable in a smooth scheme $P$,  studied by Hartshorne \cite{ha2} (and, independently, by Deligne), is defined as the hypercohomology of the completion of the De Rham complex of the formal completion of $P$ along $X$. Formal schemes play a key role in $p$-adic cohomologies (crystalline, rigid \dots) and are also algebraic models of rigid analytic spaces. These developments go back to Grothendieck with further elaborations by Raynaud, in collaboration with Bosch and L\"utkebohmert, and later work by Berthelot and de Jong. In a different vein, Strickland \cite{st} has pointed out the importance of formal schemes in the context of (stable) homotopy theory.

A particular assumption that it is almost always present in most earlier works on formal schemes is that morphisms are adic, \ie that the topology of the sheaf of rings of the initial scheme is induced by the topology of the base formal scheme. This hypothesis on a morphism of formal schemes guarantees that its fibers are usual schemes, therefore an adic morphism between formal schemes is, in the terminology of Grothendieck's school, a relative scheme over a base that is a formal scheme. But there are important examples of maps of formal schemes that do not correspond to this situation. The first example that comes into mind is the natural map $\spf(A[[X]]) \to \spf(A)$ for an adic ring $A$. This morphism has a finiteness property that had not been made explicit until \cite{AJL1} (and independently, in \cite{y}). This property is called \emph{pseudo-finite type}\footnote{In \cite{y} the terminology \emph{formally finite type} is used.}. The fact that pseudo-finite type morphisms need not be adic allows fibers that are not usual schemes, and the structure of these maps is, therefore, more complex than the structure of adic maps.
The study of smoothness and, more generally, infinitesimal lifting properties in the context of noetherian formal schemes together with this hypothesis of finiteness was embraced in general in our previous work \cite{AJP}. We should mention a preceding study of smooth morphisms under the restriction that the base is a usual scheme in \cite{y} and also the overlap of several results in \cite{AJP} and a set of results in \cite[\S 2]{LNS}, based on Nayak's 1998 thesis.

In \cite{AJP} we studied the good properties of these definitions and the agreement of their properties with the corresponding behavior for usual noetherian schemes, obtaining the corresponding statement of Zariski's Jacobian criterion for smoothness. Now we concentrate on studying properties which make sense specifically in the formal context getting information about the infinitesimal lifting properties from information present in the structure of a formal scheme. This study continues by the third author in \cite{P}  where a deformation theory for smooth morphisms is developed.

This paper can be structured roughly into three parts. The first, formed by sections 1, 2 and 3, includes preliminaries, introduces the notion of quasi-covering and the study of completion morphisms. We know of no previous reference about these matters, so we include all the needed details. They will be indispensable to state our results. The second part encompasses three sections (\ref{sec3}, \ref{sec4} and \ref{sec5}). We show that there exists  a close relationship between the infinitesimal lifting properties of  an adic  morphism and the infinitesimal lifting properties of the underlying morphism of ordinary schemes $f_{0}$. The third part (section \ref{sec6}) treats the structure theorems, which are the main results of this work. We characterize open immersions and completion morphisms in terms of the \'etale property. We classify \'etale adic coverings of a noetherian formal scheme. Finally, we give local structure theorems for unramified, \'etale and smooth maps, that show that it is possible to factor them locally into simpler maps.

Let us discuss in greater detail the contents of every section. Our framework is the category of locally noetherian formal schemes. In this category a morphism $f\colon\FX \to \FY$ can be expressed as a direct limit
\[f=\dirlim {n \in \NN} f_{n}\] 
of a family of maps of ordinary schemes using appropriate ideals of definition. The first section sets the basic notations and recalls some definitions that will be used throughout the paper. The second section deals with morphisms between locally noetherian formal schemes expressed as before as a limit
in which every map $f_{n}$ is a closed immersion of usual schemes. It is a true closed immersion of formal schemes when $f$ is adic. We treat radicial maps of formal schemes and see that the main results are completely similar to the case of usual schemes. On usual schemes, quasi-finite maps play a very important role in the understanding of the structure of \'etale maps. In the context of formal schemes there are two natural generalizations of this notion. The simplest one is \emph{pseudo-quasi-finite} (Definition \ref{defncuasifin}) --- in a few words: ``of pseudo-finite type with finite fibers". The key notion though is that of quasi-covering (Definition \ref{defncuasireves}). While both are equivalent in the context of usual schemes, the latter is a basic property of unramified and, therefore, \'etale maps between formal schemes (\emph{cf.\/} Corollaries \ref{corpnrimplcr} and \ref{corcaractlocalpe}). In section 3 we discuss  flat morphisms in the context of locally noetherian formal schemes. Next, we study morphisms of completion in this setting. They form a class of flat morphisms that are closed immersions as topological maps. Such maps will be essential for the results of the last section.

Expressing  a morphism $f\colon\FX \to \FY$ between locally noetherian formal schemes as a limit as before, it is sensible to ask about the relation that exists between the infinitesimal lifting properties of  $f$ and the infinitesimal lifting properties of the underlying morphisms of usual schemes $\{f_{n}\}_{n \in \NN}$. This is one of the main themes of the next three sections. The case of unramified morphisms is simple: $f$ is unramified if and only if $f_n$ are unramified $\forall n \in \NN$ (Proposition \ref{nrfn}). Another characterization is that $f$ is unramified if and only if $f_0$ is \emph{and} the fibers of $f$ and of $f_0$ agree (Corollary \ref{corf0imfpnr}). A consequence of this result is a useful characterization of pseudo-closed immersion as those unramified morphisms such that $f_0$ is a closed immersion (Corollary \ref{pecigf0ecnr}). Smooth morphisms are somewhat more difficult to characterize. An \emph{adic} morphism $f$ is smooth if and only if $f_0$ is and $f$ is flat (Corollary \ref{flf0l}). For a non adic morphism, one \emph{cannot} expect that the maps $f_n$ are going to be smooth when $f$ is smooth as it is shown by example \ref{exf0lisonofliso}. On the positive side, there is a nice characterization of smooth closed subschemes (Proposition \ref{ecppl}). Also, the matrix jacobian criterion holds for formal schemes, see Corollary \ref{criteriojacobiano} for a precise statement. In section \ref{sec5} we combine these results to obtain properties of \'etale morphisms. It is noteworthy to point out that a smooth pseudo-quasi-finite map need not be \'etale (Example \ref{pcf+plnope}).

The last section contains our main results. First we recover in our framework the classical fact for usual schemes \cite[(17.9.1)]{EGA44} that an open immersion is a map that is \'etale and radicial (Theorem \ref{caractencab}). We also characterize completion morphisms as those pseudo-closed immersions that are flat. This and other characterizations are given in Proposition \ref{caracmorfcompl}. Writing a locally noetherian formal scheme $\FY$ as
\[\FY= \dirlim {n \in \NN} Y_{n}\]
with respect to an ideal of definition, Proposition \ref{teorequivet} says that there is an equivalence of categories between  \'etale adic $\FY$-formal schemes  and \'etale $Y_{0}$-schemes. A special case already appears in \cite[Proposition 2.4]{y}. In fact, this result is a reinterpretation of \cite[(18.1.2)]{EGA44}. The factorization theorems are based on Theorem \ref{tppalnr} that says that an unramified morphism can be factored locally into a pseudo-closed immersion followed by an \'etale adic map. As consequences we obtain Theorem \ref{tppalet} and Theorem \ref{tppall}. They state that every smooth morphism and every \'etale  morphism factor locally as a completion morphism followed by a smooth adic morphism and an \'etale adic morphism, respectively. These results explain the local structure of smooth and \'etale morphisms of formal schemes. It has been remarked by Lipman, Nayak and Sastry in \cite[p. 132]{LNS} that this observation may simplify some developments related to Cousin complexes and duality on formal schemes.

\section{Preliminaries}\label{sec1}

We denote by $\sfn$ the category of locally noetherian formal schemes and by $\sfna$ the subcategory of locally noetherian affine formal schemes. We write $\sch$ for the category of ordinary schemes.

We assume that the reader is familiar with the basic theory of formal schemes  as is explained in \cite[\S 10]{EGA1}:  formal spectrum, ideal of definition of a formal scheme, fiber product of formal schemes, functor $M \leadsto M^\tr$ for modules over adic rings, completion of a usual scheme along a closed subscheme, adic morphisms, separated morphisms, etc.

From now on and, except otherwise indicated, every formal scheme will belong to $\sfn$. Every ring under consideration will be assumed to be noetherian. So, every complete ring and every complete module will be separated under the corresponding adic topology.

\begin{parraf}  \label{lim}
Henceforth, the following notation \cite[\S 10.6]{EGA1} will be used:
\begin{enumerate}
\item
Given $\FX \in \sfn$ and  $\CJ \subset \CO_{\FX}$ an ideal of definition  for each $n\in \NN$ we put $X_{n}:=(\FX,\CO_{\FX}/\CJ^{n+1})$ and we indicate that $\FX$ is the direct limit of the schemes $X_{n}$ by 
\[\FX=\dirlim {n \in \NN} X_{n}.\]
The ringed spaces $\FX$ and $X_{n}$ have the same underlying topological space, so we will not distinguish between a point in $\FX$ or $X_{n}$.
\item
If $f\colon\FX \to \FY$ is in $\sfn$, $\CJ \subset \CO_{\FX}$ and $\CK \subset \CO_{\FY}$ are ideals of definition such that $f^{*}(\CK)\CO_{\FX} \subset \CJ$ and $f_{n}\colon X_{n}:=(\FX,\CO_{\FX}/\CJ^{n+1}) \to Y_{n}:=(\FY,\CO_{\FY}/\CK^{n+1})$ is the morphism induced by $f$, for each $n \in \NN$, then $f$ is expressed as 
\[f = \dirlim {n\in \NN} f_{n}.\]
\item
Furthermore, given $f\colon\FX \to \FY$ a morphism in $\sfn$ and $\CK \subset \CO_{\FY}$ an ideal of definition, there exist $\CJ \subset \CO_{\FX}$ an ideal of definition such that $f^*(\CK)\CO_{\FX} \subset \CJ $. Such a pair of ideals of definition will be called $f$-\emph{compatible}. 
\end{enumerate}
\end{parraf}

\begin{parraf}
Let $f \colon \FX \to \FY$ be a morphism in $\sfn$ and let $\CJ \subset \CO_{\FX}$ and $\CK \subset \CO_{\FY}$ be $f$-compatible ideals of definition. The morphism $f$ is of \emph{pseudo-finite type (pseudo-finite)} \cite[p.7]{AJL1} if $f_{0}$ (and in fact any $f_{n}$) is of finite type (finite, respectively). Moreover, if $f$ is adic we say that $f$ is of \emph{finite type (finite)} \cite[10.13.3]{EGA1} (\cite[(4.8.2)]{EGA31}, respectively). Note that these definitions do not depend on the choice of ideals of definition.
\end{parraf}

\begin{parraf} 
\cite[Definition 2.1 and Definition 2.6]{AJP} A morphism $f\colon\FX \to \FY$ in $\sfn$ is \emph{smooth  (unramified, \'etale)} if it is of pseudo-finite type and satisfies the following lifting condition:

\begin{quote}
For all affine $\FY$-schemes $Z$ and for each closed subscheme $T\inc Z$ given by a square zero ideal $\CI \subset \CO_{Z}$ the induced map
\begin{equation*} 
\Hom_{\FY}(Z,\FX) \lto \Hom_{\FY}(T,\FX)
\end{equation*}
is surjective (injective or bijective, respectively).
\end{quote}
Moreover, if $f$ is in addition adic we say that $f$ is \emph{smooth adic (unramified adic or  \'etale adic, respectively)}.

We say that $f$ is \emph{smooth (unramified or \'etale) at $x$} if there  exists an open subset $\FU \subset \FX$ with $x \in \FU$ such that $f|_\FU$ is smooth (unramified or \'etale, respectively). 
It holds that $f$ is smooth (unramified or \'etale) if and only if 
$f$ is smooth (unramified or \'etale, respectively) at $x, \forall x \in \FX$ (\emph{cf.} \cite[Proposition 4.3, 4.1]{AJP}).
\end{parraf}

\begin{parraf}
(\emph{cf.} \cite[\S3]{AJP}) Given $f \colon \FX \to \FY$ in $\sfn$ the  \emph{differential pair of  $\FX$ over $\FY$}, $(\om^{1}_{\FX/\FY}, \hd_{\FX/\FY})$, is locally  given  by
\( (\om^{1}_{A/B},\hd_{A/B})
\)
for all open sets $\FU=\spf(A) \subset \FX$ and $\FV=\spf(B) \subset \FY$ with $f(\FU) \subset \FV$. 
The   $\CO_{\FX}$-Module $\om^{1}_{\FX/\FY}$ is called the \emph{module of  $1$-differentials of  $\FX$ over $\FY$} and the continuous $\FY$-derivation $\hd_{\FX/\FY}$ is called the \emph{canonical derivation of  $\FX$ over $\FY$}.
\end{parraf}

\begin{parraf}\label{ec}
\cite[p. 442]{EGA1} A morphism $f \colon \FZ \to \FX$ in $\sfn$ is a \emph{closed immersion} if it factors as $\FZ \xto{g} \FX' \overset{j}\inc \FX$
where $g$ is an isomorphism of $\FZ$ into a closed subscheme $\FX' \inc \FX$ of the formal scheme $\FX$ (\cite[(10.14.2)]{EGA1}). Recall from \cite[(4.8.10)]{EGA31} that a morphism $f \colon \FZ \to \FX$ in $\sfn$ is a closed immersion if it is adic and, given $\CK \subset \CO_{\FX}$ an ideal of definition of $\FX$ and $\CJ =f^{*}(\CK)\CO_{\FZ}$ the corresponding ideal of definition of $\FZ$, the induced  morphism $f_{0} \colon Z_{0} \to X_{0}$ is a closed immersion, equivalently, the induced  morphisms $f_{n} \colon Z_{n} \to X_{n}$ are closed immersions for all $n \in \NN$.

A morphism $f: \FZ \to \FX$ in $\sfn$ is an \emph{open immersion} if it factors as $\FZ \xto{g} \FX' \inc \FX$
where $g$ is an isomorphism of  $\FZ$ into an open subscheme $\FX' \inc \FX$.
\end{parraf}


\begin{defn}
Let $\FX$ be in $\sfn$, let $\CJ \subset \CO_{\FX}$ be an ideal of  definition and $x \in \FX$. We define the  \emph{topological dimension  of   $\FX$ at $x$} as
\[
\dimtop_{x}\FX = \dim_{x}X_{0}.
\]
It is easy to see that the definition does not depend on the chosen  ideal of  definition of  $\FX$. We define the \emph{topological dimension  of  $\FX$} as
\[\dimtop\FX = \sup_{x \in \FX} \dimtop_{x} \FX = \sup_{x \in \FX} \dim_{x}X_{0} =  \dim X_{0}.\]
\end{defn}

Given $A$ an $I$-adic noetherian ring, put $X = \spec(A)$ and $\FX= \spf(A)$, then $ \dimtop \FX = \dim A / I $. While the only ``visible part'' of  $\FX$ in $X = \spec(A)$ is $V(I)$, it happens that $X \setminus V(I)$ has a deep effect on the behavior of $\FX$ as we will see along this work. So apart from the topological dimension of $\FX$, it is necessary  to consider another notion of dimension that expresses part of the ``hidden" information: the algebraic dimension.

\begin{defn}
Let  $\FX$ be in $\sfn$ and let $x \in \FX$. We define the \emph{algebraic dimension of  $\FX$ at $x$} as
\[
\dim_{x}\FX = \dim \CO_{\FX,x}. 
\]
The  \emph{algebraic dimension of  $\FX$} is \[\dim\FX = \sup_{x \in \FX} \dim_{x} \FX.\]
\end{defn}

\begin{propo} \label{diaf}
If $\FX = \spf(A)$ with $A$ an $I$-adic noetherian ring then
\(\dim \FX = \dim A.\)
\end{propo}
\begin{proof}
For each  $x \in \FX$, if $\fp_{x}$ is the corresponding open prime ideal  in $A$ we have that \(\dim_{x} \FX =  \dim A_{\{\fp_{x}\}} = \dim A_{\fp_{x}}\) since $ A_{\fp_{x}} \inc A_{\{\fp_{x}\}} $ is a flat extension of local rings with the same residue field (\emph{cf.} \cite[(24.D)]{ma1}). Since $I \subset A$ is in the Jacobson radical, it holds that $\dim A = \sup_{x \in \FX} \dim A_{\fp_{x}}$, from which it follows the equality.
\end{proof}

\begin{ex} \label{exdimafdis}
Given $A$ an   $I$-adic noetherian ring and $\mathbf{T}= T_{1},\, T_{2},\, \ldots,\, T_{r}$ a finite number of indeterminates, the \emph{affine formal space of dimension $r$ over $A$} is $ \BA_{\spf(A)}^{r}=\spf(A\{\mathbf{T}\} )$ and the \emph{formal disc of dimension $r$ over $A$} is $ \BD_{\spf(A)}^{r}=\spf(A[[\mathbf{T}]])$ (see \cite[Example 1.6]{AJP}). It holds that
\[
\begin{array}{ccccc}
\dimtop \BA_{\spf(A)}^{r} &= &\dim \BA_{\spec(A/I)}^{r} &= &\dim A / I + r \\
\dimtop \BD_{\spf(A)}^{r} &= &\dim \spec(A/I)           &= &\dim A / I 
\end{array}
\]
and 
\[
\begin{array}{ccccccc}
 \dim \BA_{\spf(A)}^{r} & \underset{\textrm{\ref{diaf}}} = &\dim A\{\mathbf{T}\}  &=& \dim A + r  &\underset{\textrm{\ref{diaf}}}  = &\dim \spf(A) + r\\
 \dim \BD_{\spf(A)}^{r} &\underset{\textrm{\ref{diaf}}} = & \dim A[[\mathbf{T}]]  &=& \dim A + r  &\underset{\textrm{\ref{diaf}}}  = &\dim \spf(A) + r.
\end{array}
\]
\end{ex}

From these examples, we see that the algebraic dimension of  a formal scheme does not measure  the  dimension of the underlying topological space. 
In general, for $\FX$ in $\sfn$, $\dim_{x} \FX \ge \dimtop_{x} \FX$, for any $x \in \FX$ and, therefore
\[\dim \FX \ge \dimtop \FX.\] 
Moreover, if $\FX= \spf(A)$ with  $A$ an  $I$-adic  ring then $\dim \FX \ge \dimtop \FX + \htt (I)$.

\begin{defn}
Let $f:\FX \to \FY$ be in $\sfn$ and $y \in \FY$. The \emph{fiber of  $f$ at the  point $y$} is the formal scheme
\[f^{-1} (y) = \FX \times_{\FY} \spec(k(y)).\]
For example, if $f:\FX = \spf(B) \to \FY= \spf(A)$ is in $\sfna$ we have that $f^{-1} (y) = \spf(B \tc_{A} k(y))$.
\end{defn}

\begin{ex}
Let $\FY= \spf(A)$ be in $\sfna$ and let $\mathbf{T}= T_{1},\, T_{2},\, \ldots,\, T_{r}$ be a set of indeterminates.
If $p: \BA_{\FY}^{r} \to \FY$ is the canonical projection of the  affine formal $r$-space over $\FY$, then for all $y \in \FY$ we have that 
\[
p^{-1}(y) = \spf(A\{\mathbf{T}\}\tc_{A}k(y)) = \spec(k(y)[\mathbf{T}]) = \BA_{\spec(k(y))}^{r}.
\]
If $q: \BD_{\FY}^{r}\to\FY$ is  the canonical projection of the formal $r$-disc over $\FY$, given $y \in \FY$, it holds that 
\[
q^{-1}(y) = \spf(A[[\mathbf{T}]]\tc_{A}k(y)) = \spf(k(y)[[\mathbf{T}]]) = \BD_{\spec(k(y))}^{r}.
\]

\end{ex}
\begin{parraf} \label{fibra}
Let $f:\FX \to \FY$ be in $\sfn$ and let us consider $\CJ \subset \CO_{\FX}$ and $\CK\subset \CO_{\FY}$ $f$-compatible ideals of definition. According to \ref{lim}, 
\[f = \dirlim  {n \in \NN} (f_{n}: X_{n} \to Y_{n}).\]
Then, by \cite[(10.7.4)]{EGA1}, it holds that
\[f^{-1} (y) = \dirlim  {n \in \NN}f_{n}^{-1} (y)\] 
where $f_{n}^{-1} (y) = X_{n} \times_{Y_{n}} \spec(k(y))$, for each $n \in \NN$.

If $f$ is adic, by base-change  (\emph{cf.} \cite[1.3]{AJP}) we deduce that $f^{-1}(y) \to \spec(k(y))$ is adic so, $f^{-1} (y)$ is an ordinary scheme and $f^{-1} (y) =f_{n}^{-1} (y)$, for all $n \in \NN$.
\end{parraf}

\begin{parraf} 
We establish the following convention. Let $f:\FX \to \FY$ be in $\sfn$, $x \in \FX$ and $y= f(x)$ and assume that  $\CJ \subset \CO_{\FX}$ and $\CK\subset \CO_{\FY}$ are $f$-compatible ideals of  definition. From now  and, except otherwise indicated, whenever we consider the rings $\CO_{\FX,x}$ and $\CO_{\FY,y}$ we will associate them the  $\CJ \CO_{\FX,x}$ and $ \CK \CO_{\FY,y}$-adic topologies, respectively. And we will denote by $\widehat{\CO_{\FX,x}}$  and $\widehat{\CO_{\FY,y}}$ the completion of  $\CO_{\FX,x}$ and  $\CO_{\FY,y}$ with respect to the  $\CJ \CO_{\FX,x}$ and $ \CK \CO_{\FY,y}$-adic topologies, respectively. Note that these topologies do not depend on the choice of ideals of definition of $\FX$ and $\FY$.
\end{parraf}

\begin{defn}
Let $f:\FX \to \FY$ be in $\sfn$. Given $x \in \FX$ and $y= f(x)$, we define the  \emph{relative algebraic dimension of  $f$ at $x$} as 
\[\dim_{x} f = \dim_{x} f^{-1}(y).\]
If $\CJ \subset \CO_{\FX}$ and $\CK\subset \CO_{\FY}$ are $f$-compatible ideals of  definition, then 
\[\dim_{x} f= \dim \CO_{f^{-1}(y),x} = \dim \CO_{\FX,x} \otimes_{\CO_{\FY,y}} k(y) = \dim \widehat{\CO_{\FX,x} }\otimes_{\widehat{\CO_{\FY,y} }} k(y).\]
If the topology in $\widehat{\CO_{\FX,x} }\otimes_{\widehat{\CO_{\FY,y} }} k(y)$ is the $\CJ \widehat{\CO_{\FX,x}}$-adic then $\widehat{\CO_{\FX,x} }\tc _{\widehat{\CO_{\FY,y} }} k(y) = \widehat{\CO_{\FX,x} }\otimes_{\widehat{\CO_{\FY,y} }} k(y)$.
\end{defn}

\begin{parraf} \label{exdimalg} 
Given an adic morphism $f:\FX \to \FY$ in $\sfn$ and $f$-compatible ideals of definition $\CJ \subset \CO_{\FX}$ and $\CK\subset \CO_{\FY}$, then $\dim_{x} f  = \dim _{x} f_{0}$ for every $x \in \FX$. For example: 
\begin{enumerate}
\item \label{exdimalg1}
If $p: \BA_{\FY}^{r}:=\BA_{\ZZ}^{r} \times_{\ZZ} \FY\to \FY$ is the canonical projection of the affine formal  $r$-space over $\FY$, given $x \in \BA_{\FY}^{r}$ we have that
\[
\dim_{x} p = \dim k(y)[\mathbf{T}] = r,
\]
where $y=p(x)$.
In contrast, if  $q: \BD_{\FY}^{r}:=\BD_{\ZZ}^{r} \times_{\ZZ} \FY\to\FY$ is the canonical projection of the formal  $r$-disc over $\FY$,  $x \in \BD_{\FY}^{r}$ and $y=q(x)$ it holds  that
\[
\dim_{x} q = \dim k(y)[[\mathbf{T}]] \underset{\textrm{\ref{exdimafdis}}} = r >  \dim k(y) = 0.
\]
\item \label{exdimalg2}
If $X$ is a usual noetherian scheme  and $X'$ is a closed  subscheme  of  $X$, recall that the  morphism of  completion of  $X$ along  $X' $, $\kappa: X_{/X' } \to X$ (\cite[(10.8.5)]{EGA1})  is not adic, in general.  Note however that
\[
\dim_{x} \kappa = \dim k(x)  = 0
\]
for all  $x \in X_{/X' }$.
\end{enumerate}
\end{parraf}

\section{Pseudo-closed immersions and quasi-coverings}

\begin{defn}
 A morphism $f:\FX \to \FY$ in $\sfn$  is a  \emph{pseudo-closed immersion} if there exists $\CJ \subset \CO_{\FX}$ and $\CK \subset \CO_{\FY}$ $f$-compatible ideals of definition such that the induced morphisms of schemes $\{f_{n}:X_{n} \to Y_{n}\}_{n \in \NN}$ are closed immersions. 
 
 Note that if $f:\FX \to \FY$ is a pseudo-closed immersion, $f(\FX)$ is a closed subset of $\FY$.
 
Let us show that this definition does not depend on the choice of ideals of  definition. Being a local question, we can assume that $f: \FX= \spf(A) \to \FY = \spf(B)$ is in $\sfna$ and that $ \CJ = J^{\tr},\, \CK= K^{\tr}$ for ideals of definition $J \subset A$ and $K \subset B$ such that $KA \subset J$. Then, given another pair of ideals of definition $J' \subset A$ and $K' \subset B$ such that $\CJ' = J'^{\tr}\subset \CO_{\FX},\, \CK' = K'^{\tr}\subset \CO_{\FY}$ satisfying that $f^{*}(\CK')\CO_{\FX} \subset \CJ'$, there exists $n_{0} > 0$ such that  $J^{n_{0}} \subset J'$ and $K^{n_{0}} \subset K'$. The  morphism $B \to A$ induces the following commutative diagrams
\begin{diagram}[height=2em,w=2em,labelstyle=\scriptstyle]
B/K^{n_{0}(n+1)}	 &   \rTto    &    A/J^{n_{0}(n+1)}\\
\dTto            &            &    \dTto\\
B/K'^{n+1}       &	\rTto    &    A/J'^{n+1}\\      
\end{diagram}
and it follows that  $B/K'^{n+1}\to A/J'^{n+1}$ is surjective, for all $n \in \NN$.  Then, using \ref{ec}, it follows that the morphism $ (\FX, \CO_{\FX}/\CJ'^{n+1}) \to (\FY,\CO_{\FY}/\CK'^{n+1})$ is a  closed immersion, for all $n \in \NN$.
\end{defn}

\begin{ex}
Let $X$ be a noetherian scheme and let $X' \subset X$ be a closed subscheme defined by an  ideal $\CI \subset \CO_{X}$. The  morphism of  completion $X_{/X'} \xto{\kappa} X$ of  $X$ along  $X'$ (\cite[(10.8.5)]{EGA1}) is expressed as
\[
\dirlim {n \in \NN} \left((X' ,\CO_{X}/\CI^{n+1}) \xto{ \kappa_{n}} (X,\CO_{X})\right),
\]
therefore, it is a pseudo-closed immersion.
\end{ex}

Notice that an adic pseudo-closed immersion is a closed immersion (\emph{cf.} \ref{ec}). However, to be a pseudo-closed immersion is not a topological property:
\begin{ex}
 Given $K$ a field, let $p: \BD_{\spec(K)}^{1}  \to \spec(K)$ be the canonical projection. If we consider the ideal of  definition $\langle T\rangle^{\tr}$, of  $\BD_{\spec(K)}^{1}$ then $p_{0} = 1_{\spec(K)}$ is a  closed immersion. However,  the morphisms \[p_{n}: \spec(K[T] / \langle T \rangle  ^{n+1}) \to \spec(K)\] are not closed immersions, for all $ n > 0$ and, thus, $p$ is not a pseudo-closed immersion.
\end{ex}

\begin{propo} 
Let $f:\FX \to \FY$ and $g:\FY \to \FS$ be two morphisms in $\sfn$. It holds that:
\begin{enumerate}
\item 
If $f$ and $g$ are pseudo-closed immersions then $g \circ f$ is a pseudo-closed immersion.

\item 
If $f$ is a pseudo-closed immersion, given $h:\FY' \to \FY$ in $\sfn$  we have that $\FX_{\FY'} = \FX \times_{\FY} \FY'$ is in $\sfn$ and  that $f':\FX_{\FY'} \to \FY'$ is a pseudo-closed immersion.

\end{enumerate}
\end{propo}

\begin{proof} 
As for (1) let $\CJ \subset \CO_{\FX}$, $\CK \subset \CO_{\FY}$ and $\CL \subset \CO_{\FS}$ be ideals of  definition such that $\CJ$ and $\CK$ are $f$-compatible, $\CK$ and $\CL$ are $g$-compatible and consider the corresponding expressions for $f$ and $g$ as direct limits:
\[
f =\dirlim {n \in \NN} (X_{n} \xto{f_{n}} Y_{n}) \qquad g= \dirlim {n \in \NN} (Y_{n} \xto{g_{n}}S_{n})
\]
Since
\[g \circ f = \dirlim {n \in \NN}  g_{n} \circ f _{n}\] the assertion follows from the stability under composition  of  closed immersions in $\sch$. 
Let us show (2).  Take $\CK' \subset \CO_{\FY'}$ an ideal of  definition with $h^{*}(\CK)\CO_{\FY'} \subset \CK'$ and such that, by \ref{lim},  
\[h= \dirlim {n \in \NN} (h_{n}: Y'_{n} \to Y_{n}).\] Then by \cite[(10.7.4)]{EGA1} we have that 
\[
\begin{matrix}
\begin{diagram}[height=2.5em,w=2.5em,labelstyle=\scriptstyle]
\FX_{\FY'}	& \rTto^{f'}	& \FY' &   \\
\dTto       &         &\dTto^h   & \\
\FX	        &\rTto^f  &	\FY &	\\
\end{diagram}
& = \qquad & \dirlim {n \in \NN}&
\left(
\begin{diagram}[height=2.5em,w=2.5em,labelstyle=\scriptstyle]
 X_{n} \times_{Y_{n}} Y'_{n} & \rTto^{f'_{n}}&Y'_{n}\\
  \dTto		&   &\dTto^{h_{n}}\\
  X_{n}	&		\rTto^{f_{n}} & Y_{n}\\
\end{diagram}
\right)
\end{matrix}
\]
By hypothesis, $f_{n}$ is a  closed immersion and since  closed immersions in $\sch$ are stable under base-change we have that $f'_{n}$ is a  closed immersion of  noetherian schemes, $\forall n\in \NN$. Finally, since $f$ is a morphism of pseudo-finite type, from \cite[Proposition 1.8.(2)]{AJP} we have that $\FX_{\FY'}$ is in $\sfn$.
\end{proof}


Next we turn to the study of radicial morphisms in the context of formal schemes. This notion will allow us later (Theorem \ref{caractencab}) to give a characterization of open immersions in terms of \'etale morphisms.

\begin{defn}\label{rad}
A morphism $f:\FX \to \FY$ in $\sfn$ is \emph{radicial} if given $\CJ \subset \CO_{\FX}$ and $\CK \subset \CO_{\FY}$ $f$-compatible ideals of  definition the  induced morphism of  schemes $f_{0}:X_{0} \to Y_{0}$ is radicial. 

Given $x \in \FX$, the residue fields of the local rings $\CO_{\FX,x}$ and  $\CO_{X_{0},x}$ agree and analogously for $\CO_{\FY,f(x)}$ and  $\CO_{Y_{0},f_0(x)}$. Therefore the definition of radicial morphisms does not depend on the chosen ideals of definition of $\FX$ and $\FY$.
\end{defn}
           
\begin{parraf} \label{prad}
From the sorites of radicial morphisms in $\sch$  it follows that:
\begin{enumerate}
\item \label{prad1}
Radicial morphisms are stable under composition  and noetherian base-change.
\item \label{prad2}
Every monomorphism is radicial. So, open immersions, closed immersions and pseudo-closed immersions are radicial morphisms. 
\end{enumerate}
\end{parraf}

The notion of quasi-finite morphism of usual schemes (see \cite[Definition (6.11.3)]{EGA1}) is based on the equivalence between several conditions for morphisms between schemes (Corollaire (6.11.2) in \emph{loc.~cit.}) that are no longer equivalent in the full context of formal schemes. Specifically, we study two notions that generalize that of quasi-finite morphism of usual schemes. They will play a basic role in understanding the structure of unramified and \'etale morphisms in $\sfn$.

\begin{defn} \label{defncuasifin}
Let $f\colon\FX \to \FY$ be a pseudo-finite type morphism in $\sfn$. We say that $f$ is \emph{pseudo-quasi-finite} if there exist $\CJ \subset \CO_{\FX}$ and $\CK \subset \CO_{\FY}$ $f$-compatible ideals of definition such that $f_{0}$ is quasi-finite. And $f$ is \emph{pseudo-quasi-finite at $x \in \FX$} if there exists an open neighborhood $x \in\FU \subset \FX$ such that $f|_{\FU}$ is pseudo-quasi-finite.

Notice that if $f:\FX \to \FY$ is a pseudo-quasi-finite morphism (in $\sfn$) then,  \emph{for all $f$-compatible ideals of  definition} $\CJ \subset \CO_{\FX}$ and $\CK \subset \CO_{\FY}$, the  induced morphism of  schemes $f_{0}:X_{0} \to Y_{0}$ is quasi-finite.

As an immediate consequence  of  the analogous properties in $\sch$ we have that:

\begin{enumerate}
\item
The underlying sets of the fibers of a pseudo-quasi-finite morphism are finite.
\item \label{soritpcf1}
Closed immersions, pseudo-closed immersions and open immersions are pseudo-quasi-finite.
\item \label{soritpcf1bis}
Pseudo-finite and finite morphisms are pseudo-quasi-finite.
\item \label{soritpcf2}
If $f:\FX \to \FY$ and $g:\FY \to \FS$ are pseudo-quasi-finite morphisms, then so is $g \circ f$.
\item  \label{soritpcf3}
If $f:\FX \to \FY$ is pseudo-quasi-finite, given $h:\FY' \to \FY$ a morphism  in $\sfn$ we have that $f':\FX_{\FY'} \to \FY'$ is pseudo-quasi-finite.
\end{enumerate}
\end{defn}

In $\sch$ it is the case that a morphism is \'etale if and only if it is smooth and quasi-finite. However, we will show that in $\sfn$ not  every smooth and pseudo-quasi-finite morphism is \'etale (see Example \ref{pcf+plnope}). That is why we introduce a stronger notion than pseudo-quasi-finite morphism and that also generalizes quasi-finite morphisms in $\sch$: the quasi-coverings.

\begin{defn} \label{defncuasireves}
Let $f:\FX \to \FY$ be a pseudo-finite type morphism in $\sfn$. The  morphism $f$ is a \emph{quasi-covering} if  $\CO_{\FX,x} \tc_{\CO_{\FY,f(x)}} k(f(x))$ is a finite type $k(f(x))$-module, for all $x \in \FX$. We say that $f$ is a \emph{quasi-covering at $x \in \FX$} if there exists an open $\FU \subset \FX$ with $x \in \FU$ such that  $f|_{\FU}$ is a  quasi-covering.

We reserve the word \emph{covering} for a dominant (\ie with dense image) quasi-covering. These kind of maps will play no role in the present work but they are important, for instance, in the study of finite group actions on formal schemes.
\end{defn}

\begin{ex}
If $X$ is a locally noetherian scheme and $X' \subset X $ is a closed subscheme  the  morphism of  completion $\kappa: \FX=X_{/X'} \to X$ is a quasi-covering. In fact, for all $x \in \FX$ we have that 
\[\CO_{\FX,x} \tc_{\CO_{X,\kappa(x)}} k(\kappa(x)) = k(\kappa(x)).\]
\end{ex}

\begin{lem} \label{soritcr}
We have the following: 
\begin{enumerate}
\item  \label{soritcr1}
Closed immersions,  pseudo-closed immersions and open immersions are  quasi-coverings.
\item  \label{soritcr2}
If $f:\FX \to \FY$ and $g:\FY \to \FS$  are quasi-coverings, the  morphism $g \circ f$ is a quasi-covering.
\item   \label{soritcr3}
If $f:\FX \to \FY$ is a quasi-covering, and $h:\FY' \to \FY$ a morphism  in $\sfn$, then $f':\FX_{\FY'} \to \FY'$ is a quasi-covering.
\end{enumerate}
\end{lem}

\begin{proof}
Immediate.
\end{proof}
\begin{propo} \label{cuairevdim0}
If $f:\FX \to \FY$ is a quasi-covering  in $x \in \FX$ then:
\[
\dim_{x} f =0
\]
\end{propo}
\begin{proof}
It is a consequence of the fact that $\CO_{\FX,x} \tc_{\CO_{\FY,f(x)}} k(f(x))$ is an artinian ring.
\end{proof}

\begin{rem}
Observe that given $\CJ \subset \CO_{\FX}$ and $\CK \subset \CO_{\FY}$ ideals of  definition such that $f^{*}(\CJ) \CO_{\FX}\subset \CK$, for all $x \in \FX$ it holds that 
\[\CO_{\FX,x} \tc_{\CO_{\FY,f(x)}} k(f(x)) = \invlim {n \in \NN} \CO_{X_{n},x} \otimes_{\CO_{Y_{n},f_n(x)}} k(f(x)).\]
Over usual schemes quasi-coverings and pseudo-quasi-finite morphisms are equivalent notions. More generally we have the following.
\end{rem}

\begin{propo}
Let $f:\FX \to \FY$ be a morphism in $\sfn$. If $f$ is a quasi-covering, then  it is pseudo-quasi-finite. Furthermore, if  $f$ is adic then the converse holds.
\end{propo}

\begin{proof}
Suppose that $f$ is a quasi-covering and let $\CJ \subset \CO_{\FX}$ and $\CK\subset \CO_{\FY}$ be $f$-compatible ideals of definition. Given $x \in \FX$ and $y=f(x)$, $\CO_{\FX,x} \tc_{\CO_{\FY,y}} k(y)$ is a finite $k(y)$-module  and, therefore, 
\[\frac{\CO_{X_{0},x}}{\fm_{Y_{0},y}\CO_{X_{0},x}} = \frac{\CO_{\FX,x}}{ \CJ \CO_{\FX,x}} \otimes_{\CO_{Y_{0},y}} k(y)\]
is $k(y)$-finite, so it follows that $f$ is pseudo-quasi-finite.
 
If  $f$ is an adic  morphism, $f^{-1}(y)=f_{0}^{-1}(y)$ for each $y \in \FY$, so 
\[\CO_{X_{0},x}/\fm_{Y_{0},y}\CO_{X_{0},x} = \CO_{\FX,x} \tc_{\CO_{\FY,y}} k(y)\]
for all $x \in \FX$ with  $y=f(x)$. If $f$ is moreover pseudo-quasi-finite, it follows from \cite[Corollaire (6.11.2)]{EGA1} that $f$ is a quasi-covering.
\end{proof}

\begin{cor}
Every finite  morphism $f:\FX \to \FY$  in $\sfn$ is a quasi-covering.
\end{cor}

\begin{proof}
Finite morphisms are adic and pseudo-quasi-finite. Therefore the  result  is  consequence of the last proposition.
\end{proof}

Nevertheless, by the next example, not every pseudo-finite morphism is a quasi-covering and, therefore, \emph{pseudo-quasi-finite} does not imply \emph{quasi-covering} for morphisms in $\sfn$.
\begin{ex} 
For $r > 0$, the canonical projection $p: \BD_{\FX}^{r} \to \FX$  is not  a quasi-covering since 
\[\dim_{x} p \underset{\textrm{\ref{exdimalg}.(\ref{exdimalg1})}} = r > 0 \quad
\forall  x \in \FX.\]
But  considering an appropriate pair of ideals of definition, the scheme map $p_{0}= 1_{X_{0}}$ is  finite.
\end{ex}

\begin{parraf}
In short, we have  the  following diagram of  strict implications (with the conditions that imply adic morphism in italics):
\[
\begin{matrix}
\textrm{\emph{closed immersion}}&\imp &\textrm{\emph{finite}}&\imp&
\textrm{quasi-covering}\\
\Downarrow &    &\Downarrow & &\Downarrow\\
\textrm{pseudo-closed immersion}&\imp &\textrm{pseudo-finite}&\imp & \textrm{pseudo-quasi-finite}\\
\end{matrix}
\]
\end{parraf}

\section{Flat morphisms and completion morphisms} \label{sec2}

In  the first part  of  this section we discuss  flat morphisms in $\sfn$. Whenever a  morphism \[f = \dirlim {n \in \NN} f_{n}\] is adic, the  local criterion of flatness for  formal schemes (Proposition \ref{clp}) relates the flat character  of  $f$ and that of  the  morphisms $f_{n}$, for all $n \in \NN$.  In absence of  the adic hypothesis this relation does not hold, though (Example \ref{excompl1}). In the second part, we study  the morphisms of  completion in $\sfn$, a  class of flat morphisms that are pseudo-closed immersions (so, they are closed immersions as topological maps). Even though the construction of the completion of a formal scheme along a closed formal subscheme is a natural one, it has not been systematically developed in the basic references about formal schemes. Morphisms of completion will be an essential ingredient in  the main theorems of Section  \ref{sec6}, namely, Theorems \ref{tppalnr}, \ref{tppalet} and \ref{tppall}.

\begin{parraf}\label{caracterizlocalplanos}
A morphism $f:\FX \to \FY$ is \emph{flat at $x \in \FX$} if $\CO_{\FX,x}$ is a flat $\CO_{\FY,f(x)}$-module. We say that \emph{$f$ is flat} if it is flat at $x$, for all  $x \in \FX$.

Given $\CJ \subset \CO_{\FX}$ and $\CK\subset \CO_{\FY}$ $f$-compatible ideals of  definition, by \cite[III, \S5.4, Proposition 4]{b} the following are equivalent:
\begin{enumerate}
\item
$f$ is flat at $x \in \FX$.
\item
$\CO_{\FX,x}$ is a flat $\CO_{\FY,f(x)}$-module.
\item
$\widehat{\CO_{\FX,x}}$ is a flat $\CO_{\FY,f(x)}$-module.
\item
$\widehat{\CO_{\FX,x}}$ is a flat $\widehat{\CO_{\FY,f(x)}}$-module.
\end{enumerate} 
\end{parraf}

\begin{ex} \label{excompl1}
Let $K$ be a field let $\BA_{K}^{1}= \spec(K[T])$ and consider the  closed subset $X' =V( \langle T \rangle) \subset  \BA_{K}^{1}$. The  canonical morphism  of  completion of   $\BA_{K}^{1}$ along  $X'$
\[\BD_{K}^{1}  \xto{\kappa} \BA_{K}^{1}\] 
is flat but,  the morphisms 
\[\spec(K[T]/ \langle T \rangle ^{n+1}) \xto{\kappa_{n}} \BA_{K}^{1}\]
are not flat, for every $n \in \NN$.
\end{ex}

\begin{propo} \label{clp}
(Local flatness criterion for formal schemes.)
Given an \emph{adic} morphism $f\colon\FX \to \FY$ in $\sfn$, and an ideal of definition $\CK \subset \CO_{\FY}$, then $\CJ =  f^{*}(\CK)\CO_{\FX} \subset \CO_{\FX}$ is an ideal of definition. Let $\{f_{n}\colon X_{n} \to Y_{n}\}_{n \in \NN}$ be the morphisms  induced  by $f$ through $\CK$ and $\CJ$. The following assertions are equivalent:
\begin{enumerate}
\item
The morphism $f$ is flat.
\item
The morphism $f_{n}$ is flat, for all $n \in \NN$.
\item
The morphism $f_{0}$ is flat.
\end{enumerate}
\end{propo}

\begin{proof}
We may suppose that $f: \FX=\spf(A) \to \FY=\spf(B)$ is in $\sfn$. Then if $\CK=K^{\tr}$ for an ideal of definition $K \subset B$, we have that $\CJ= (KA)^{\tr}$  and the proposition is a consequence of \cite[Lemma 7.1.1]{AJL1} and of the local flatness criterion for rings (\emph{cf.} \cite[Theorem 22.3]{ma2}).
\end{proof}


Associated to  a (usual) locally noetherian scheme $X$  and  a closed subscheme of  $X' \subset X$  there is a locally noetherian formal scheme $X_{/X'}$, called completion of  $X$ along  $X'$ and,  a canonical morphism $\kappa: X_{/X'} \to X$ (\cite[(10.8.3) and (10.8.5)]{EGA1}). Next, we define the completion of  a formal scheme  $\FX$ along a closed formal subscheme $\FX'\subset \FX$.
  
\begin{defn} \label{defcom}
Let $\FX$ be in $\sfn$ and let $\FX' \subset \FX$ be a closed formal subscheme defined by a coherent ideal $\CI$ of  $\CO_{\FX}$. Given an ideal of definition $\CJ$ of  $\FX$  we define the completion of a sheaf $\CF$ on $\FX$ over $\FX'$, denoted by $\CF_{/ \FX'}$, as the restriction to $\FX'$ of the sheaf
\[
\invlim {n \in \NN} \frac{\CF}{(\CJ+\CI)^{n+1}\CF}.
\] 
The definition does not depend neither on the chosen  ideal of  definition $\CJ$ of  $\FX$ nor on the coherent ideal $\CI$ that defines $\FX'$.

We define the \emph{completion of  $\FX$ along  $\FX'$}, and it will be denoted $\FX_{/ \FX'}$, as the  topological ringed space whose underlying topological space is $\FX'$ and whose sheaf of topological rings is $\CO_{\FX_{/\FX'}}$. 

It is easy to check that $\FX_{/ \FX'}$ satisfies the hypothesis of  \cite[(10.6.3) and (10.6.4)]{EGA1},  from which we deduce that:
\begin{enumerate}
\item The formal scheme
 $\FX_{/ \FX'}$ is locally noetherian.
 \item The ideal
 $(\CI+\CJ)_{/\FX'} \subset \CO_{\FX_{/\FX'}}$ defined by the restriction to $\FX'$ of the sheaf 
 \[\invlim {n \in \NN} \frac{\CJ+\CI}{(\CJ+\CI)^{n+1}}\] 
 is an ideal of  definition of  $\FX_{/\FX'}$.
 \item It holds that
 $\CO_{\FX_{/\FX'}}/ ((\CI+\CJ)_{/\FX'})^{n+1}$ agrees with the restriction to $\FX'$ of the sheaf $\CO_{\FX}/ (\CJ+\CI)^{n+1}$ for every $n \in \NN$.
\end{enumerate}
\end{defn}

\begin{parraf} \label{limcomple}
With the above notations, if $Z_{n} = (\FX',  \CO_{\FX}/ (\CJ+\CI)^{n+1})$  for all $n \in \NN$, by \ref{lim} we have that \[\FX_{/ \FX'} = \dirlim {n \in \NN}Z_{n}\]
For each $n \in \NN$, let $X_{n} =(\FX, \CO_{\FX}/\CJ^{n+1})$ and $X'_{n} =(\FX', \CO_{\FX}/(\CJ^{n+1}+\CI))$. The canonical  morphisms 
\[\frac{\CO_{\FX}}{\CJ^{n+1}} \epi \frac{\CO_{\FX}}{(\CJ+\CI)^{n+1}} \epi \frac{\CO_{\FX}}{\CJ^{n+1} + \CI}\] 
provide the  closed immersions of  schemes $X'_{n} \xto{j_{n}} Z_{n} \xto{\kappa_{n}} X_{n}$, such that the diagram, whose vertical maps are the obvious closed immersions,
\begin{diagram}[height=2em,w=2em,labelstyle=\scriptstyle]
X'_{m} & \rTto^{j_{m}} & Z_{m}	& \rTto^{\kappa_{m}} & X_{m} \\
\uTto  &               & \uTto	&                    & \uTto\\
X'_{n} & \rTto^{j_{n}} & Z_{n}	& \rTto^{\kappa_{n}} & X_{n}\\
\end{diagram}
commutes, for all $m \ge n \ge 0$.
Then by \ref{lim} we have the canonical morphisms in $\sfn$
\[
\FX' \xto{j} \FX_{/ \FX'} \xto{\kappa} \FX 
\]
where $j$ is a  closed immersion (see  \ref{ec}). The  morphism $\kappa$ as topological map is the inclusion and it is called  \emph{morphism of  completion of  $\FX$ along $\FX'$}. 
\end{parraf}

\begin{rem}
Observe that $\kappa$ is adic only if $\CI$ is contained in a ideal of  definition of  $\FX$, in which case $\FX= \FX_{/ \FX'}$ and $\kappa=1_{\FX}$.
\end{rem}

\begin{parraf} \label{compl}
If $\FX= \spf(A)$ is in $\sfna$ with $A$ a  $J$-adic noetherian ring, and $\FX'= \spf(A/I)$ is a closed formal subscheme  of  $\FX$, then 
\[\ga(\FX_{/ \FX'},\CO_{\FX_{/ \FX'}}) =  \invlim {n \in \NN} \frac{A}{(J+I)^{n+1}}=: \HA\]
and from \cite[(10.2.2) and (10.4.6)]{EGA1} we have that 
$\FX_{/ \FX'}=\spf(\HA)$ and the morphisms $\FX' \xto{j} \FX_{/ \FX'} \xto{\kappa} \FX$ correspond to the natural continuous  morphisms $A \to \HA \to A /I $. 
\end{parraf}

\begin{propo} \label{caractcom} 
Given $\FX$ in $\sfn$ and $\FX'$ a closed formal subscheme of  $\FX$, the  morphism of  completion  $\kappa: \FX_{/ \FX'} \to \FX$ is a pseudo-closed immersion and \'etale (and therefore, from \cite[Proposition 4.8]{AJP}, it is flat).
\end{propo}
\begin{proof}
With  the notations of  \ref{limcomple} we have that 
\[\kappa= \dirlim {n \in \NN} \kappa_{n}.\]
Since $\kappa_{n}$ is a  closed immersion for all $n \in \NN$, it follows that $\kappa$ is a pseudo-closed immersion.
In order to prove that $\kappa$ is an \'etale morphism we may suppose that
$\FX=\spf(A)$ and  $\FX'= \spf(A/I)$, where $A$ is a $J$-adic noetherian ring. Note that $\FX_{/ \FX'}=\spf(\HA)$ where $\HA$ is the  completion  of  $A$ for the $(J+I)$-adic topology and, therefore, is \'etale over $A$. By \cite[2.2]{AJP}, $\kappa$ is an \'etale morphism. 
\end{proof}

\begin{rem}
In Theorem  \ref{caracmorfcompl} we will see that the converse holds: every  flat  pseudo-closed immersion is a morphism of  completion.
\end{rem}

\begin{parraf}  \label{complmorf}
Given $f:\FX \to \FY$ in $\sfn$, let $\FX' \subset \FX$ and  $\FY' \subset \FY$ be closed formal subschemes given by ideals $\CI \subset \CO_{\FX}$ and $\CL \subset \CO_{\FY}$ such that $f^{*}(\CL)\CO_{\FX} \subset \CI$, that is, $f(\FX')\subset \FY'$. If $\CJ \subset \CO_{\FX}$ and $\CK \subset \CO_{\FY}$ are $f$-compatible ideals of  definition, let us denote for all $n \in \NN$
\begin{align*}
X_{n}  &= (\FX,  \CO_{\FX}/\CJ^{n+1}),       & \quad 
Y_{n}  &= (\FY,  \CO_{\FY}/\CK^{n+1}),         \\
Z_{n}  &= (\FX', \CO_{\FX}/(\CJ+\CI)^{n+1}), & \quad 
W_{n}  &= (\FY', \CO_{\FY}/(\CK+\CL)^{n+1})    \\
X'_{n} &= (\FX', \CO_{\FX}/(\CJ^{n+1}+\CI))  & \quad
Y'_{n} &= (\FY', \CO_{\FX}/(\CK^{n+1}+\CL)).
\end{align*}
Then the  morphism $f$ induces the following commutative  diagram of  locally noetherian schemes where the oblique maps are the canonical immersions:
\begin{diagram}[height=1.5em,w=2em,labelstyle=\scriptstyle]
X_{n}&&\rTto^{f_{n}}&&Y_{n}&&&\\
&\rdTto&&&\uTto&\rdTto&&\\
\uTto^{\kappa_{n}}&&X_{m}&\rTto^{f_{m}}&\HonV&&Y_{m}&\\
&&\uTto^{\kappa_{m}}&&\vLine^{\kappa'_{n}}&&&\\
Z_{n}&\hLine&\VonH&\rTto^{\hf_{n}}&W_{n}&&\uTto_{\kappa'_{m}}&\\
&\rdTto&&&\uTto&\rdTto&&\\
\uTto^{j_{n}}&&Z_{m}&\rTto^{\hf_{m}}&\HonV&&W_{m} & \\
&&\uTto^{j_{m}}&&\vLine^{j'_{n}}&&&\\
X'_{n}&\hLine&\VonH&\rTto^{f'_{n}}& Y'_{n}&&\uTto_{j'_{m}}&\\
&\rdTto&&&&\rdTto&&\\
&&X'_{m}&&\rTto^{f'_{m}}&& Y'_{m}& \\
\end{diagram}
for all $m \ge n \ge 0$. 
Note that
\[
f' = \dirlim {n \in \NN} f'_{n}
\]
is the restriction $f|_{\FX'} \colon \FX' \to \FY'$.
Applying the direct limit over ${n \in \NN}$ we obtain  a morphism 
\[\hf : \FX_{/ \FX'} \to \FY_{/ \FY'}\]
 in $\sfn$, such that  the  following diagram is commutative:
\begin{equation} \label{diagrcom}
\begin{diagram}[height=2em,w=2em,labelstyle=\scriptstyle]
\FX            & \rTto^{f}         & \FY \\
\uTto^{\kappa} &                   & \uTto^{\kappa'} \\
\FX_{/ \FX'}   & \rTto^{\hf}       & \FY_{/ \FY'} \\
\uTto          &                   & \uTto \\
\FX'           & \rTto^{f|_{\FX'}} & \FY'
\end{diagram}
\end{equation}
We will call $\hf$ the \emph{completion of  $f$ along  $\FX'$ and $\FY'$}.
\end{parraf}

\begin{parraf} \label{complafin}
Suppose that $f: \FX=\spf(A) \to \FY=\spf(B)$ is in $\sfna$ and that $\FX'= \spf(A/I)$ and $\FY'= \spf(B/L)$ with $LA \subset I$. If $J \subset A$ and $K \subset B$ are ideals of  definition such that $KA \subset J$, the  morphism $\hf : \FX_{/ \FX'} \to \FY_{/ \FY'}$ corresponds  to the   morphism induced by $B \to A$
\[\widehat{B}\to \HA\]
(\emph{cf.} \cite[(10.4.6)]{EGA1}) where $\HA$ is the completion of  $A$ for  the $(I+J)$-adic topology and $\widehat{B}$ denotes the  completion of  $B$ for the  $(K+L)$-adic topology.
\end{parraf}

\begin{propo} \label{cccom}
Given $f:\FX \to \FY$ in $\sfn$, let $\FY' \subset \FY$ be a closed formal subscheme and $\FX' = f^{-1}(\FY')$. Then, 
\[\FX_{/ \FX'} = \FY_{/ \FY'} \times_{\FY} \FX.\]
\end{propo}

\begin{proof}
We may restrict to the case in which $\FX=\spf(A)$, $\FY=\spf(B)$ and $\FY'= \spf(B/L)$ are affine formal schemes and $J \subset A$ and $K \subset B$ are ideals of  definition such that $KA \subset J$. By hypothesis,  $\FX'= \spf(A/LA)$, so $\FX_{/ \FX'} = \spf(\HA)$ where $\HA$ is the  completion of  $A$ for the  $(J+LA)$-adic topology. On the other hand, $\FY_{/ \FY'} = \spf(\HB)$ where $\HB$ denotes the completion of  $B$ for the $(K+L)$-adic  topology and  it holds that 
\[\HB \tc_{B} A= B\tc_{B} A=\HA,\]
since $J+(K+L)A = J+KA+LA= J+LA$, the result follows.
\end{proof}

\begin{propo} \label{complsorit} 
Given $f:\FX \to \FY$ in $\sfn$, let us consider closed formal subschemes $\FX' \subset \FX$ and $\FY' \subset \FY$ such that  $f(\FX')\subset \FY'$.
\begin{enumerate}
\item
Let $\CP$ be one of  the following properties of  morphisms in $\sfn$: 
\begin{center}
\emph{pseudo-finite type,  pseudo-finite, pseudo-closed immersion, pseudo-quasi-finite, quasi-covering, flat, separated, radicial, smooth, unramified, \'etale.}
\end{center}
If  $f$ satisfies $\CP$, then so does $\hf$. 
\item
Moreover, if $\FX' = f^{-1}(\FY')$, let $\CQ$ be one of  the following properties of  morphisms in $\sfn$: 
\begin{center}
\emph{adic, finite type,  finite,  closed immersion, smooth adic, unramified adic, \'etale adic.}
\end{center}
Then, if $f$ satisfies $\CQ$, then so does $\hf$. 
\end{enumerate}
\end{propo}
\begin{proof}

Suppose that $f$ is flat and let us prove that $\hf$ is flat. The question is local so we may assume   $f: \FX=\spf(A) \to \FY=\spf(B)$  in $\sfna$, $\FX'= \spf(A/I)$ and $\FY'= \spf(B/L)$ with $LA \subset I$. Let $J \subset A$ and $ K \subset B$ be ideals of  definition such that $KA \subset J$ and, $\HA$ and $\widehat{B}$  the completions of $A$ and $B$ for the topologies given by $(I+J)\subset A$ and $(K+L) \subset B$, respectively. By \cite[III, \S5.4, Proposition 4]{b} we have that the  morphism $\widehat{B}\to \HA$ is  flat and, from  \ref{complafin} and  \cite[Lemma 7.1.1]{AJL1} it follows that $\hf$ is flat.

Suppose that $f$ satisfies any of the other properties $\CP$ and let us prove that $\hf$ inherits them using the commutativity of the diagram
\begin{diagram}[height=2em,w=2em,labelstyle=\scriptstyle]
\FX             & \rTto^{f}   & \FY\\
\uTto^{ \kappa} &             & \uTto^{ \kappa'}\\
\FX_{/ \FX'}    & \rTto^{\hf} & \FY_{/ \FY'} \\
\end{diagram}
where the vertical arrows are morphisms of completion. Since all of these properties $\CP$ are stable under composition and a morphism of completion satisfies $\CP$ (Proposition \ref{caractcom}) we have that  $\CP$ holds for $f \circ \kappa = \kappa' \circ \hf$.  
If $\CP$ is smooth, unramified or \'etale the result is immediate from \cite[Proposition 2.13]{AJP}. 
 
If $\CP$ is any of the other properties, then closed immersions verify $\CP$ and $\CP$ is stable under composition and under base-change in $\sfn$. Therefore, since $\kappa' \circ \hf$ has $\CP$ and $\kappa'$ is separated (Proposition \ref{caractcom}), by the analogous argument in $\sfn$ to the one in $\sch$  \cite[(5.2.7), $i),\, ii) \Rightarrow iii)$]{EGA1} we get that $\hf$ also satisfies $\CP$.

Finally, if $f$ is adic, from Proposition \ref{cccom} and from \cite[1.3]{AJP}, we deduce that $\hf$ is adic. Then, if $\CQ$ is any of the properties in statement (2) and $f$ satisfies $\CQ$, by (1) so does $\hf$.
\end{proof}


\section{Unramified morphisms} \label{sec3}
Let $f:\FX \to \FY$ be a morphism of locally noetherian formal schemes. Given $\CJ \subset \CO_{\FX}$ and $\CK \subset \CO_{\FY}$ $f$-compatible ideals of definition, express $f$ as a limit
\[f:\FX \to \FY= \dirlim {n \in \NN} (f_{n}:X_{n} \to Y_{n}).\]
We begin relating the  unramified character of $f:\FX \to \FY$ 
and that of the underlying ordinary scheme morphisms  $\{f_{n}\}_{n \in \NN}$. 

\begin{propo}\label{nrfn}
With the previous notations, the  morphism $f$ is unramified  if and only if  $f_{n}:X_{n} \to Y_{n}$  is unramified, for all $n \in \NN$.
\end{propo}

\begin{proof}
Notice that both conditions in the statement imply that $f$ is a pseudo-finite type morphism. Applying  \cite[Proposition 4.6]{AJP} we have  to show that $\om^{1}_{\FX/\FY} = 0$ is equivalent to $\Omega^{1}_{X_{n}/Y_{n}}=0$, for all $n \in \NN$. If $\om^{1}_{\FX/\FY} =0$,  by the Second Fundamental Exact Sequence (\cite[Proposition 3.13]{AJP}) for the morphisms \[X_{n} \inc \FX \xto{f} \FY,\] we have that $\Omega^{1}_{X_{n}/ \FY}=0$, for all $n \in \NN$.  From  the First Fundamental Exact Sequence (\cite[Proposition 3.10]{AJP}) associated to the morphisms \[X_{n} \xto{f_{n}} Y_{n} \inc \FY,\] it follows that $\Omega^{1}_{X_{n}/Y_{n}}=0$. The converse follows from the identification
\[\om^{1}_{\FX/\FY} =\invlim {n \in \NN} \Omega ^{1}_{X_{n}/Y_{n}}  \]
(\emph{cf}. \cite[\S 1.9]{AJP} ).
\end{proof}

\begin{cor} \label{pecimplnoram}
With the previous notations, if the morphisms $f_{n}:X_{n} \to Y_{n}$ are  immersions  for all $n \in \NN$, then $f$ is unramified. 
\end{cor}

In the class of  adic  morphisms in $\sfn$ the  following proposition provides a criterion, stronger than the last result,  to  determine when a morphism $f$ is unramified.

\begin{propo} \label{fnrf0nr}
Let $f:\FX \to \FY$ be an \emph{adic} morphism in $\sfn$ and let $\CK \subset \CO_{\FY}$ be an ideal of  definition. Write \[f=\dirlim {n \in \NN} f_{n}\] by taking ideals of  definition $\CK \subset \CO_{\FY}$ and $\CJ=f^{*}(\CK)\CO_{\FX} \subset \CO_{\FX}$. The  morphism $f$ is unramified if and only if the induced morphism $f_{0}:X_{0} \to Y_{0}$ is unramified.
\end{propo}

\begin{proof}
If $f$ is unramified by Proposition  \ref{nrfn} we have that $f_{0}$ is unramified. 
Conversely, suppose that  $f_{0}$ is unramified and let us prove that $\om^{1}_{\FX/\FY}=0$. The question is local so we may assume that $f:\FX=\spf(A) \to \FY= \spf(B)$ is in $\sfna$ and that $\CJ=J^{\tr},\, $ with $J \subset A$ an ideal of  definition. By hypothesis $\Omega^{1}_{X_{0}/Y_{0}}=0$  and thus, since $f$ is adic it holds that 
\begin{equation} \label{eqmodifceromodiff0cero}
\om^{1}_{\FX/\FY} \otimes_{\CO_{\FX}} \CO_{X_{0}}\underset{\textrm{\cite[3.8]{AJP}}} = \Omega^{1}_{X_{0}/Y_{0}}=0.
\end{equation}
Then by the equivalence  of categories \cite[(10.10.2)]{EGA1}, the last equality says that $\om^{1}_{A/B}/J \om^{1}_{A/B}=0$. Since $A$ is a $J$-adic ring it holds that  $J$ is contained in the Jacobson  radical  of  $A$. Moreover, \cite[Proposition 3.3]{AJP} implies that $\om^{1}_{A/B}$ is a  finite type $A$-module. From Nakayama's  lemma we deduce that $\om^{1}_{A/B}=0$ and therefore, $\om^{1}_{\FX/\FY}= (\om^{1}_{A/B})^{\tr}= 0$. Applying  \cite[Proposition 4.6]{AJP} it follows that $f$ is unramified.
\end{proof}

The following example illustrates that in the  non adic  case  the  analogous of  the last  proposition does not hold.

\begin{ex} \label{framfonram}
Let  $K$ be a field and $p:\BD^{1}_{K} \to \spec(K)$ be the projection morphism of  the formal disc of  dimension $1$ over $\spec(K)$. By \cite[Example 3.14]{AJP} we have that $\om^{1}_{p} =(K[[T]] \hd T)^{\tr}$ and therefore, $\BD^{1}_{K}$ is  ramified over $K$ (\cite[Proposition 4.6]{AJP}).  However, given the  ideal of  definition $\langle T \rangle \subset K[[T]]$ the induced morphism $p_{0} = 1_{\spec(K)}$ is unramified.
\end{ex}

Let us consider for a morphism $f\colon\FX \to \FY$ in $\sfn$ the notation established at the beginning of the section. In view of the example, our next goal will be to determine when the morphism $f$ such that  $f_{0}$  is unramified but $f$ itself is not necessarily adic, is unramified (Corollary \ref{corf0imfpnr}). In order to do that, we will need some results that describe the local behavior of unramified morphisms. Next, we provide local characterizations of  unramified  morphisms in $\sfn$, generalizing the  analogous properties in the category of  schemes (\emph{cf.} \cite[(17.4.1)]{EGA44}).

\begin{propo} \label{caraclocalpnr}
Let $f:\FX \to \FY$ be a morphism in $\sfn$ of  pseudo-finite type. For $x \in \FX$ and $y=f(x)$ the following conditions are equivalent:
\begin{enumerate}
\item[(1)]
 $f$ is unramified  at $x$.
\item[(2)]
$f^{-1}(y)$ is an unramified  $k(y)$-formal scheme  at $x$.
\item[(3)]
$\fm_{\FX,x}\widehat{\CO_{\FX,x}} = \fm_{\FY,y}\widehat{\CO_{\FX,x}}$ and $k(x)|k(y)$ is a finite  separable extension.
\item[(4)]
$\om^{1}_{\CO_{\FX,x}/ \CO_{\FY,y}}  =0$.
\item[$(4')$]
$(\om^{1}_{\FX/\FY})_{x}=0$.
\item[(5)]
$\CO_{\FX,x}$ is a formally unramified $\CO_{\FY,y}$-algebra  for  the  adic topologies.
\item[$(5')$]
$\widehat{\CO_{\FX,x}}$ is a formally unramified $\widehat{\CO_{\FY,y}}$-algebra   for  the  adic topologies.
\end{enumerate}
\end{propo}

\begin{proof}
Keep the notation from the beginning of this section and write
\[f:\FX \to \FY= \dirlim {n \in \NN} (f_{n}:X_{n} \to Y_{n})
\]
(1) $\dimp$ (2) By Proposition \ref{nrfn}, $f$ is unramified  at $x$ if and only if all the morphisms $f_{n}:X_{n} \to Y_{n}$ are unramified at $x$. Applying  \cite[(17.4.1)]{EGA44}, this is equivalent to $f_{n}^{-1}(y)$ being an unramified $k(y)$-scheme at $x$, for all $n \in \NN$, which  is also equivalent   to  \[f^{-1} (y) \underset{\textrm{\ref{fibra}}}= \dirlim {n \in \NN} f_{n}^{-1}(y)\] being an unramified $k(y)$-formal scheme  at $x$.

(1) $\imp$ (3) The assertion (1) is equivalent to $f_{n}:X_{n} \to Y_{n}$ being unramified at $x$, for all $n \in \NN$, and from \cite[\emph{loc.~cit.}]{EGA44} it follows that  $k(x)|k(y)$  is a finite separable extension, and that $\fm_{X_{n},x} = \fm_{Y_{n},y} \CO_{X_{n},x}$, for all $n \in \NN$. Hence, 
\[
\fm_{\FX,x} \widehat{\CO_{\FX,x}} = \invlim {n \in \NN} \fm_{X_{n},x} =  \invlim {n \in \NN} \fm_{Y_{n},y} \CO_{X_{n},x} = \fm_{\FY,y}\widehat{\CO_{\FX,x}}.
\]

(4) $\dimp$ ($4'$) By \cite[Proposition 3.3]{AJP} it holds that $(\om^{1}_{\FX/\FY})_{x}$ is a finite type $\CO_{\FX,x}$-module   and therefore,  
\[
\om^{1}_{ \CO_{\FX,x}/ \CO_{\FY,y}} = \widehat{(\om^{1}_{\FX/\FY})_{x}} = (\om^{1}_{\FX/\FY})_{x} \otimes_{\CO_{\FX,x}} \widehat{\CO_{\FX,x}}.
\]
Then, since $\widehat{\CO_{\FX,x}}$ is a faithfully flat $\CO_{\FX,x}$-algebra, $\om^{1}_{ \CO_{\FX,x}/ \CO_{\FY,y}}=0$ if and only if $(\om^{1}_{\FX/\FY})_{x}=0$.

(3) $\imp$ (4) Since $k(x)|k(y)$  is a finite separable extension  we have that $\Omega^{1}_{k(x)/k(y)}=0$ and from \cite[Proposition 3.3]{AJP}  $ \om^{1}_{\CO_{\FX,x}/ \CO_{\FY,y}}= \widehat{(\om^{1}_{\FX/\FY})_{x}}$ is a  finite type $\widehat{\CO_{\FX,x}}$-module. Therefore, it holds  that \[\om^{1}_{\CO_{\FX,x}/ \CO_{\FY,y}} \otimes_{\widehat{\CO_{\FX,x}}} k(x)= \om^{1}_{(\CO_{\FX,x} \otimes_{\CO_{\FY,y}} k(y))/k(y)} =\Omega^{1}_{k(x)/k(y)}=0.\] 
By   Nakayama's lemma, $\om^{1}_{ \CO_{\FX,x}/ \CO_{\FY,y}} =0$.

(4) $\dimp$ (5) It is straightforward from \cite[(\textbf{0}, 20.7.4)]{EGA41}.

(5) $\dimp$ ($5'$) Immediate.

($4'$) $\imp$ (1) Since $\om^{1}_{\FX/\FY} \in \coh(\FX)$ (\cite[Proposition 3.3]{AJP}), assertion $(4')$ implies that there exists  an open subset $\FU \subset \FX$  with $x \in \FU$ such that  $(\om^{1}_{\FX/\FY})|_{\FU}=0$ and therefore,  by \cite[Proposition 4.6]{AJP} we have that $f$ is unramified  at $x$.
\end{proof}

\begin{cor} \label{corcaraclocalpnr}
Let $f:\FX \to \FY$ be a pseudo-finite type morphism in $\sfn$. The following conditions are equivalent:
\begin{enumerate}
\item[(1)]
 $f$ is unramified.
\item[(2)]
For all $x \in \FX$, $f^{-1}(f(x))$ is an unramified $k(f(x))$-formal scheme  at $x$.
\item [(3)]
For all $x \in \FX$, $\fm_{\FX,x}\widehat{\CO_{\FX,x}} = \fm_{\FY,f(x)}\widehat{\CO_{\FX,x}}$ and $k(x)|k(f(x))$ is a finite separable extension. 
\item[(4)]
$\om^{1}_{\CO_{\FX,x}/ \CO_{\FY,f(x)}}  =0$, for all $x \in \FX$. 
\item[$(4')$]
For all $x \in \FX$, $(\om^{1}_{\FX/\FY})_{x} =0$.
\item[(5)]
For all $x \in \FX$, $\CO_{\FX,x}$ is a formally unramified $\CO_{\FY,f(x)}$-algebra for the  adic topologies.
\item[$(5')$]
For all $x \in \FX$, $\widehat{\CO_{\FX,x}}$ is a formally unramified  $\widehat{\CO_{\FY,f(x)}}$-algebra for the  adic topologies.
\end{enumerate}

\end{cor}

\begin{cor} \label{corpnrimplcr}
Let $f:\FX \to \FY$ be a pseudo-finite type morphism in $\sfn$. If $f$ is unramified  at $x \in \FX$, then $f$ is a quasi-covering at $x$.
\end{cor}

\begin{proof}
By assertion (3) of  Proposition \ref{caraclocalpnr} we have that \[\CO_{\FX,x} \tc_{\CO_{\FY,f(x)}} k(f(x))=k(x)\] with $k(x)|k(f(x))$ a finite extension and therefore, $f$ is a quasi-covering at $x$ (see Definition \ref{defncuasireves}).
\end{proof}

\begin{cor} \label{pnrdim0} 
Let $f:\FX \to \FY$ be a pseudo-finite type morphism in $\sfn$. If $f$ is unramified  at $x \in \FX$, then $\dim_{x} f =0$.
\end{cor}

\begin{proof}
It is straightforward from the previous Corollary and Proposition \ref{cuairevdim0}.
\end{proof}

\begin{propo} \label{fonrydislocal}
Let $f:\FX \to \FY$ be a pseudo-finite type morphism in $\sfn$. Given $x \in \FX$ and $y=f(x)$ the following conditions are equivalent:
\begin{enumerate}
\item
$f$ is unramified  at $x$.
\item
$f_{0}:X_{0} \to Y_{0}$ is unramified at $x$ and $\widehat{\CO_{\FX,x}} \otimes_{\widehat{\CO_{\FY,y}}} k(y) = k(x)$.
 \end{enumerate}
\end{propo}

\begin{proof}
If $f$ is unramified at $x$, then $f_{0}$ is unramified at $x$ (Proposition \ref{nrfn}). Moreover,  assertion (3) of  Proposition \ref{caraclocalpnr} implies that $\widehat{\CO_{\FX,x}} \otimes_{\widehat{\CO_{\FY,y}}} k(y) = k(x)$ so (1) $\Rightarrow$ (2) holds. Let us prove  that (2) $\Rightarrow$ (1). Since $f_{0}$  is unramified at $x$  we have that $k(x)|k(y)$ is a finite separable extension (\emph{cf.} \cite[(17.4.1)]{EGA44}). From  the equality $\widehat{\CO_{\FX,x}} \otimes_{\widehat{\CO_{\FY,y}}} k(y) = k(x)$ we deduce that $\fm_{\FX,x}\widehat{\CO_{\FX,x}} = \fm_{\FY,y}\widehat{\CO_{\FX,x}}$. Thus,  the  morphism $f$  and the  point $x$ satisfy assertion (3) of Proposition  \ref{caraclocalpnr} and it follows that $f$ is unramified  at $x$.
\end{proof}

Now we are ready to state the non adic version of Proposition \ref{fnrf0nr}:

\begin{cor}  \label{corf0imfpnr}
Given $f:\FX \to \FY$ a morphism in $\sfn$ of  pseudo-finite type let $\CJ \subset \CO_{\FX}$ and $\CK \subset \CO_{\FY}$ be $f$-compatible ideals of  definition and let $f_0 \colon X_0 \to Y_0$ be the induced morphism. The following conditions are equivalent:
\begin{enumerate}
\item
The morphism $f$ is unramified.
\item
The morphism $f_{0}$ is unramified and, for all $x \in \FX$, $f^{-1}(y)=f_{0}^{-1}(y)$ with $y=f(x)$.
\end{enumerate}
\end{cor}

\begin{proof}
Suppose that $f$ is unramified and fix $x \in \FX$ and $y=f(x)$. By Proposition  \ref{fonrydislocal} we have that $f_{0}$ is unramified and that  $\widehat{\CO_{\FX,x}} \otimes_{\widehat{\CO_{\FY,y}}} k(y) = k(x)$. Therefore, $\CJ(\widehat{\CO_{\FX,x}} \otimes_{\widehat{\CO_{\FY,y}}} k(y) )= 0$ and applying  Lemma \ref{fibradiscr} we deduce that $f^{-1}(y)=f_{0}^{-1}(y)$. Conversely,  suppose that  (2) holds and let us show that given $x \in \FX$, the  morphism $f$ is unramified  at $x$. If $y=f(x)$, we have that $f_{0}^{-1}(y)$ is an unramified   $k(y)$-scheme at $x$ (\emph{cf.} \cite[(17.4.1)]{EGA44}) and since $f^{-1}(y)=f_{0}^{-1}(y)$, from   Proposition \ref{caraclocalpnr} it follows that $f$ is unramified  at $x$.
\end{proof}

\begin{lem} \label{fibradiscr} 
Let $A$ be a $J$-adic noetherian ring such that  for all open prime ideals $\fp \subset A$, $J_{\fp}=0$. Then $J=0$ and therefore, the  $J$-adic topology  in $A$ is the discrete topology.  
\end{lem}

\begin{proof}
Since every maximal ideal $\fm \subset A$ is open for the $J$-adic topology, we have that $J_{\fm}=0$, for all maximal ideal $\fm \subset A$, so $J=0$.
\end{proof}

\begin{parraf}
As a consequence of Corollary \ref{corf0imfpnr} it holds that: 
\begin{itemize}
\item
If $f:\FX \to \FY$ is an unramified morphism in $\sfn$ then  $f^{-1}(y)$ is a usual scheme for all $x \in \FX$ where  $y=f(x)$.
\item
In  Corollary \ref{corcaraclocalpnr} assertion (2) may be written:
\item[$(2')$]
\emph{For all $x \in \FX,\, y=f(x)$, $f^{-1}(y)$ is a unramified $k(y)$-scheme  at $x$.}

\end{itemize}
\end{parraf}

From Proposition \ref{caraclocalpnr} we obtain the  following result, in which we provide a description of  pseudo-closed immersions that will be used  in the characterization of  completion morphisms (Theorem \ref{caracmorfcompl}).

\begin{cor} \label{pecigf0ecnr} 
Given $f:\FX \to \FY$ in $\sfn$, let $\CJ \subset \CO_{\FX}$ and $\CK \subset \CO_{\FY}$ be $f$-compatible ideals of definition and express \[f=\dirlim {n\in \NN} f_{n}.\] The  morphism $f$ is a pseudo-closed immersion if and only if $f$ is unramified  and $f_{0}:X_{0} \to Y_{0}$ is a  closed immersion.
\end{cor}

\begin{proof}
If $f$ is a pseudo-closed immersion, by  Corollary  \ref{pecimplnoram} it follows that $f$ is unramified. Conversely, suppose that $f$ is unramified and that $f_{0}$ is a  closed immersion and let us show that $f_{n}:X_{n} \to Y_{n}$ is a  closed immersion, for each $n \in \NN$. By \cite[(4.2.2.(ii))]{EGA1} it suffices to prove that, for all $x \in \FX$ with $y=f(x)$, the  morphism  $\CO_{Y_{n},y} \to \CO_{X_{n},x}$ is surjective, for all $n \in \NN$. Fix $x \in \FX$, $ y=f(x) \in \FY$ and $n \in \NN$. Since $f_{0}$ is a  closed immersion, by \cite[\emph{loc. cit.}]{EGA1},   we have that $\CO_{Y_{0},y} \to \CO_{X_{0},x}$ is surjective and therefore, $\spf(\widehat{\CO_{\FX, x}}) \to \spf(\widehat{\CO_{\FY, y}})$ is a  pseudo-finite morphism, so, the  morphism $\CO_{Y_{n},y} \to \CO_{X_{n},x}$ is  finite. On the other hand, the morphism $f$ is unramified therefore by Proposition \ref{nrfn} we get that $f_{n}$ is unramified and applying Proposition \ref{caraclocalpnr}   we obtain that $\fm_{Y_{n},y} \CO_{X_{n},x}= \fm_{X_{n},x}$. Then by  Nakayama's lemma we conclude that  $\CO_{Y_{n},y} \to \CO_{X_{n},x}$ is a surjective morphism.
\end{proof}

\section{Smooth morphisms} \label{sec4}

The contents of  this section can be structured in two parts. In the first part we study the relationship between the smoothness of  a morphism \[f= \dirlim {n \in \NN} f_{n}\] in $\sfn$ and the smoothness of the ordinary scheme morphisms $\{f_{n}\}_{n \in \NN}$. In the second part, we provide a local factorization for smooth morphisms (Proposition  \ref{factpl}). In this section we also prove in Corollary \ref{criteriojacobiano} the matrix  Jacobian criterion, that is a useful explicit condition in terms of a matrix rank for determining whether a closed subscheme of the affine formal space or of the affine formal disc is smooth or not.

\begin{propo} \label{lfn}
Given $f:\FX \to \FY$ in $\sfn$, let $\CJ \subset \CO_{\FX}$ and $\CK \subset \CO_{\FY}$ be $f$-compatible ideals of  definition and  write
\[f=\dirlim {n\in \NN} f_{n}.\]
If $f_{n}:X_{n} \to Y_{n}$ is smooth, for all $n \in \NN$, then $f$ is smooth.
\end{propo}

\begin{proof}
By \cite[Proposition 4.1]{AJP} we may assume that $f$ is in $\sfna$. 
Let $Z$ be an affine scheme, consider a morphism $w: Z \to \FY$, a closed $\FY$-subscheme $T \inc Z$ given by a square zero ideal and a $\FY$-morphism $u: T \to \FX$. Since $f$ and $w$ are morphisms of  affine formal schemes we find  an integer $m \ge 0$ such that  $w^{*}(\CK^{m+1}) \CO_{Z} =0$ and $u^{*}(\CJ^{m+1})\CO_{T}=0$  and therefore $u$ and $w$ factors as $T \xto{u_{m}} X_{m} \xto{i_{m}} \FX $ and $Z \xto{w_{m}} Y_{m} \xto{i_{m}} \FY $, respectively. Since $f_{m}$ is formally smooth, there exists a $Y_{m}$-morphism $v_{m}: Z \to X_{m}$ such that  the  following diagram is commutative
\begin{diagram}[height=2.5em,w=2em,labelstyle=\scriptstyle]
T             & \rTinc               &                     & Z \\
\dTto^{u_{m}} &                      & \ldTto(3,2)^{v_{m}} & \dTto^{w_{m}}\\ 
X_{m}         & \rTto^{\qquad f_{m}} &                     & Y_{m}\\ 
\dTto^{i_{m}} &                      &                     & \dTto\\ 
\FX           & \rTto^{\qquad f}     &                     & \FY\\
\end{diagram}
Thus the  $\FY$-morphism $v:=i_{m} \circ v_{m}$ satisfies that $v|_{T}=u$ and then, $f$ is formally smooth. Moreover,  since $f_{0}$ is a finite type morphism, it holds that $f$ is of  pseudo-finite type and therefore,  $f$ is smooth. 
\end{proof}

\begin{cor} \label{ladfn}
Let $f:\FX \to \FY$ be an \emph{adic}  morphism  in $\sfn$ and consider $\CK \subset \CO_{\FY}$ an  ideal of  definition.  
The morphism $f$ is smooth  if and only if all the scheme morphisms $\{f_{n}:X_{n} \to Y_{n}\}_{n \in \NN}$, determined by the ideals of  definition $\CK \subset \CO_{\FY}$ and $\CJ=f^{*}(\CK)\CO_{\FX}$, are  smooth.
\end{cor}

\begin{proof}
If $f$ is  adic, by \cite[(10.12.2)]{EGA1}, we have that for each $n \in \NN$, the diagram
\begin{diagram}[height=2em,w=2em,labelstyle=\scriptstyle]
\FX   & \rTto^{f}     & \FY\\
\uTto &               & \uTto\\ 
X_{n} & \rTto^{f_{n}} & Y_{n}\\ 
\end{diagram}
is a cartesian square.
Then by base-change (\cite[Proposition 2.9 (2)]{AJP}) we have that $f_{n}$ is smooth, for all $n \in \NN$. The  converse follows from the previous proposition.
\end{proof}
Next example shows us that the  converse of  Proposition  \ref{lfn} does not hold in general. 

\begin{ex} \label{peynopen}
Let $K$ be a field and $ \BA_{K}^{1}= \spec(K[T])$. For the  closed subset $X =V( \langle T \rangle) \subset  \BA_{K}^{1}$,   Proposition \ref{caractcom} implies that the  canonical completion morphism 
\[\BD_{K}^{1}  \xto{\kappa} \BA_{K}^{1}\] 
of $\BA_{K}^{1}$ along  $X$ is \'etale. However,  picking in $\BA_{K}^{1}$ the  ideal of  definition $0$,   the morphisms 
\[\spec(K[T]/ \langle T \rangle ^{n+1}) \xto{\kappa_{n}} \BA_{K}^{1}\]
are not flat, whence it follows that $\kappa_{n}$ can not be smooth for all $n \in \NN$ (see \cite[Proposition 4.8]{AJP}).
\end{ex}

Our next goal will be to determine the relation between smoothness of  a morphism 
\[f= \dirlim {n \in \NN} f_{n}\]  and that of  $f_{0}$ (Corollaries \ref{flf0l} and \ref{corf0imfpl}). In order to do that, we need to characterize smoothness locally.

\begin{propo} \label{pligplafibr}
Let $f:\FX \to \FY$ be a pseudo-finite type morphism in $\sfn$. Given $x \in \FX$ and $y = f(x)$ the following conditions are equivalent:
\begin{enumerate}
\item
The morphism $f$ is smooth  at $x$.
\item
$\CO_{\FX,x}$ is a  formally smooth $\CO_{\FY,y}$-algebra for the adic topologies.
\item
$\widehat{\CO_{\FX,x}}$ is a formally smooth $\widehat{\CO_{\FY,y}}$-algebra  for the adic topologies.
\item
The morphism $f$ is  flat at $x$ and $f^{-1}(y)$ is a $k(y)$-formal scheme smooth  at $x$.
\end{enumerate}
\end{propo}

\begin{proof}
The question is  local and $f$ is of  pseudo-finite type, so we may assume  that $f: \FX= \spf(A) \to \FY= \spf(B)$ is in $\sfna$, with 
$A = B\{T_{1},\, \ldots ,T_{r}\}[[Z_{1},\ldots, Z_{s}]]/I$ and 
$I \subset B':= B\{T_{1},\, \ldots ,T_{r}\}[[Z_{1},\ldots, Z_{s}]]$ 
an ideal (\cite[Proposition 1.7]{AJP}). Let $\fp \subset A$ be the open prime  ideal corresponding to $x$, let $\fq \subset B'$ be the open prime such that  $\fp = \fq/I$ and let $\fr \subset B$ be the  open prime  ideal corresponding to   $y$. 

(1) $\imp$ (3) Replacing $\FX$ by a sufficiently small  open neighborhood   of  $x$ we may suppose  that $A$ is a  formally smooth $B$-algebra.  Then, by \cite[(\textbf{0}, 19.3.5)]{EGA41} we have that $A_{\fp}$ is a formally smooth $B_{\fr}$-algebra  and    \cite[(\textbf{0}, 19.3.6)]{EGA41} implies  that $\widehat{\CO_{\FX,x}}=\widehat{A_{\fp}}$ is a  formally smooth $\widehat{\CO_{\FY,y}}=\widehat{B_{\fr}}$-algebra.

(2) $\dimp$ (3) It is a consequence of \cite[(\textbf{0}, 19.3.6)]{EGA41}.

(3) $\imp$ (1) By  \cite[(\textbf{0}, 19.3.6)]{EGA41}, assertion (3) is equivalent 
to  $A_{\fp}$ being a formally  smooth $B_{\fr}$-algebra. Then Zariski's Jacobian criterion  (\cite[Proposition 4.14]{AJP}  implies that the  morphism of  $\widehat{A_{\fp}}$-modules
\[
\widehat{\frac{I_{\fq}}{I_{\fq}^{2}}} \to \Omega^{1}_{B'_{\fq}/B_{\fr}} \tc_{B'_{\fq}} A_{\fp}
\]
is right invertible. Since $\widehat{A_{\fp}}$ is a faithfully  flat $A_{\{\fp\}}$-algebra and the $A_{\{\fp\}}$-module$(\om^{1}_{B'/B} \otimes_{B'} A)_{\{\fp\}}$ is projective  (see \cite[Proposition 4.8]{AJP}), it holds that the  morphism 
\[
\left(\frac{I}{I^{2}}\right)_{\{\fp\}} \to 
(\om^{1}_{B'/B} \otimes_{B'} A)_{\{\fp\}}
\]
is right invertible by \cite[(\textbf{0}, 19.1.14.(ii))]{EGA41}. From  the equivalence of categories \cite[(10.10.2)]{EGA1} we find an open  subset $\FU \subset \FX$ with $x \in \FU$ such that  the  morphism
\[
\left(\frac{I}{I^{2}}\right)^{\tr} \to \om^{1}_{\BD^{s}_{\BA^{r}_{\FY}}/\FY} \otimes _{\CO_{\BD^{s}_{\BA^{r}_{\FY}}}} \CO_{\FX}
\]
is right  invertible over $\FU$. Now, by Zariski's Jacobian criterion for formal schemes (\cite[Corollary 4.15]{AJP}) it follows that $f$ is smooth  in $\FU$.

(3)  $\imp$ (4)  By   \cite[(\textbf{0}, 19.3.8)]{EGA41} we have that $\widehat{\CO_{\FX,x}}$ is  a formally smooth $\widehat{\CO_{\FY,y}}$-algebra for the topologies given by the maximal ideals. Then it follows from \cite[(\textbf{0}, 19.7.1)]{EGA41} that $\widehat{\CO_{\FX,x}}$ is  $\widehat{\CO_{\FY,y}}$-flat and by \ref{caracterizlocalplanos}, $f$ is flat at $x$. Moreover from  \cite[(\textbf{0}, 19.3.5)]{EGA41} we deduce that $\widehat{\CO_{\FX,x}}\otimes_{\widehat{\CO_{\FY,y}}} k(y)$ is a formally smooth $k(y)$-algebra  for the adic topologies or, equivalently, by (3) $\dimp$ (1), $f^{-1}(y)$ is a $k(y)$-formal scheme smooth  at $x$.

(4) $\imp$ (3) 
By \ref{caracterizlocalplanos} we have that $A_{\fp}$ is a flat $B_{\fr}$-module and therefore, it holds that
\begin{equation} \label{ecuacioncita}
0 \to \frac{I_{\fq}}{\fr I_{\fq}} \to \frac{B'_{\fq}}{\fr B'_{\fq}} \to \frac{A_{\fp} }{\fr A_{\fp}} \to 0
\end{equation}
is an exact  sequence. 
On the other hand, since  $f^{-1}(y)$ is a $k(y)$-formal scheme smooth  at $x$, from (1) $\Rightarrow$ (2) we deduce that $\widehat{\CO_{\FX,x}}\otimes_{\widehat{\CO_{\FY,y}}} k(y)$ is a formally smooth $k(y)$-algebra  for the adic topologies or, equivalently by  \cite[(\textbf{0}, 19.3.6)]{EGA41}, $A_{\fp}/\fr A_{\fp}$ is a formally smooth $k(\fr)$-algebra  for the adic topologies. Applying Zariski's Jacobian criterion (\cite[Proposition 4.14]{AJP}), we have that the  morphism 
\[
\widehat{\frac{I_{\fq}}{I_{\fq}^{2}}} \otimes_{B_{\fr}} k(\fr) \to (\om^{1}_{B'/B})_{\fq} \tc_{B'_{\fq}} A_{\fp} \otimes_{B_{\fr}} k(\fr)
\]
is right  invertible. Now, since  $(\om^{1}_{B'/B})_{\fq}$ is a projective $B'_{\fq}$-module (see \cite[Proposition 4.8]{AJP})
by \cite[(\textbf{0}, 6.7.2)]{EGA1} we obtain that 
\[
\widehat{\frac{I_{\fq}}{I_{\fq}^{2}}}  \to \om^{1}_{B'_{\fq}/B_{\fr}} \tc_{\widehat{B'_{\fq}}} \widehat{A_{\fp}}
\]
is right invertible. Again, by the Zariski's Jacobian criterion, $A_{\fp}$ is a formally smooth  $B_{\fr}$-algebra for the adic topologies or, equivalently by  \cite[(\textbf{0}, 19.3.6)]{EGA41}, $\widehat{A_{\fp}}$ is a  formally smooth $\widehat{B_{\fr}}$-algebra.

\end{proof}

\begin{cor} \label{corpligplafibr}
Let $f:\FX \to \FY$ be a pseudo-finite type morphism in $\sfn$. The following conditions are equivalent:
\begin{enumerate}
\item
The morphism $f$ is smooth.
\item
For all $x \in \FX$, $\CO_{\FX,x}$ is a formally smooth $\CO_{\FY,f(x)}$-algebra for the adic topologies.
\item
For all $x \in \FX$, $\widehat{\CO_{\FX,x}}$ is a formally smooth  $\widehat{\CO_{\FY,f(x)}}$-algebra for the adic topologies.
\item
The morphism $f$ is  flat and $f^{-1}(f(x))$ is a $k(f(x))$-formal scheme smooth  at $x$, for all $x \in \FX$.
\end{enumerate}
\end{cor}

\begin{cor} \label{flf0l}
Let $f:\FX \to \FY$ be  an  \emph{adic}  morphism in $\sfn$ and let $\CK \subset \CO_{\FY}$ be an ideal of  definition. Put \[f=\dirlim {n \in \NN} f_{n}\] using  the ideals of  definition $\CK \subset \CO_{\FY}$ and $\CJ=f^{*}(\CK)\CO_{\FX} \subset \CO_{\FX}$. Then, the  morphism $f$ is smooth if and only if it is  flat and  the  morphism $f_{0}:X_{0} \to Y_{0}$ is smooth.
\end{cor}

\begin{proof}
Since $f$ is  adic,  the  diagram
\begin{diagram}[height=2em,w=2em,labelstyle=\scriptstyle]
\FX   & \rTto^{f}     & \FY \\
\uTto &               & \uTto \\ 
X_{0} & \rTto^{f_{0}} & Y_{0} \\ 
\end{diagram}
is a cartesian square (\cite[(10.12.2)]{EGA1}).
If $f$ is smooth,
by base-change it follows  that $f_{0}$ is smooth. Moreover by \cite[Proposition 4.8]{AJP} we have that $f$ is flat. Conversely, if $f$ is adic, by \ref{fibra}, we have that $f^{-1}(f(x))=f_{0}^{-1}(f(x))$, for all $x \in \FX$. Therefore, since $f_{0}$ is smooth, by base-change it holds that $f^{-1}(f(x))$ is a $k(f(x))$-scheme smooth at $x$, for all $x \in \FX$ and applying    Corollary \ref{corpligplafibr} we conclude that $f$ is smooth. 
\end{proof}

The upcoming example shows that the  last result is not true without assuming the \emph{adic} hypothesis for the morphism $f$.

\begin{ex} \label{exf0lisonofliso}
Given $K$ a field, let $\mathbb{P}^{n}_{K}$ be the $n$-dimensional projective space and $X \subset \mathbb{P}^{n}_{K}$ a closed subscheme that is not smooth  over $K$. If we denote by $(\mathbb{P}^{n}_{K})_{/ X}$ the completion of  $\mathbb{P}^{n}_{K}$ along  $X$, by Proposition  \ref{complsorit} we have that the  morphism 
\[
(\mathbb{P}^{n}_{K})_{/ X} \xto{f} \spec(K)
\]
is smooth  but $f_{0} : X \to \spec(K)$ is not smooth.
\end{ex}

\begin{cor} \label{corf0imfpl}
Given $f:\FX \to \FY$ a morphism  in $\sfn$ let $\CJ \subset \CO_{\FX}$ and $\CK \subset \CO_{\FY}$ be $f$-compatible ideals of definition. Write 
\[f=\dirlim {n \in \NN} f_{n}.\] 
If $f$ is flat, $f_{0}: X_{0} \to Y_{0}$ is a smooth morphism  and $f^{-1}(f(x))=f_{0}^{-1}(f(x))$, for all $x \in \FX$, then $f$ is smooth.
\end{cor}

\begin{proof}
Since $f_{0}$ is smooth and $f^{-1}(y)=f_{0}^{-1}(y)$ for all $y = f(x)$ with $x \in \FX$,  we deduce that $f^{-1}(y)$ is a smooth  $k(y)$-scheme. Besides, by hypothesis $f$ is flat and  Corollary \ref{corpligplafibr} implies that $f$ is smooth.
\end{proof}
Example \ref{exf0lisonofliso} illustrates that the  converse of the last  corollary does not hold.

Every smooth morphism $f: X \to Y$ in $\sch$ is locally a composition of  an \'etale morphism  $U \to \BA^{r}_{Y}$ and a projection $\BA^{r}_{Y} \to Y$. Proposition \ref{factpl} generalizes this fact for smooth morphisms in $\sfn$. The same result has already appeared stated in local form in \cite[Proposition 1.11]{y}. We include it here for completeness.

\begin{propo} \label{factpl}
Let $f:\FX \to \FY$ be a pseudo-finite type morphism in $\sfn$. The  morphism $f$ is smooth  at $x \in \FX$ if  and only if there exists an open subset $\FU \subset \FX$ with $x \in \FU$ such that  $f|_{\FU}$ factors as
\[
\FU \xto{g} \BA^{n}_{\FY} \xto{p} \FY
\]
where $g$ is \'etale, $p$ is the canonical projection and $n = \rg (\om^{1}_{\CO_{\FX,x}/\CO_{\FY,f(x)}})$.
\end{propo}

\begin{proof}
As this is a local  question, we may assume that $f: \FX= \spf(A) \to \FY= \spf(B)$ is a smooth morphism   in $\sfna$. By  \cite[Proposition 4.8]{AJP} and by \cite[(10.10.8.6)]{EGA1} we have that $\om^{1}_{A/B}$ is a projective $A$-module of  finite type and therefore, if $\fp \subset A$ is the open prime ideal corresponding to $x$, there exists $h\in A \setminus  \fp$  such that  $\ga(\fD(h), \om^{1}_{\FX/\FY}) = \om^{1}_{A_{\{h\}}/B}$ is a free $A_{\{h\}}$-module of  finite type. Put $\FU = \spf(A_{\{h\}})$. Given $\{\hd a_{1},\hd a_{2},\ldots, \hd a_{n}\}$ a basis of  $\om^{1}_{A_{\{h\}}/B}$ consider the  morphism  of  $\FY$-formal schemes 
\[
\FU \xto{g} \BA^{n}_{\FY}=\spf(B\{T_{1},T_{2},\ldots,T_{n}\} )
\] 
defined  by the continuous morphism of topological $B$-algebras 
\[
\begin{array}{ccc}
B\{T_{1},T_{2},\ldots,T_{n}\} & \to & A_{\{h\}}\\
T_{i} &\leadsto  &a_{i}
\end{array}
\]
See \cite[(10.2.2) and (10.4.6)]{EGA1}.
The morphism $g$ satisfies  that $f|_{\FU} = p \circ g$. Moreover, we deduce that $g^{*} \om^{1}_{\BA^{n}_{\FY}/\FY}\cong \om^{1}_{\FX/\FY}$ (see the definition of  $g$) and by  \cite[Corollary 4.13]{AJP} we have that $g$ is \'etale. 
\end{proof}

\begin{cor} \label{dimrango}
Let $f:\FX \to \FY$ be a smooth morphism   at $x \in \FX$ and $y = f(x)$. Then
\[
\dim_{x} f = \rg (\om^{1}_{\CO_{\FX,x}/\CO_{\FY,y}}).
\]
\end{cor}

\begin{proof}
Put  $n= \rg (\om^{1}_{\CO_{\FX,x}/\CO_{\FY,y}})$. By Proposition  \ref{factpl} there exists $\FU \subset \FX$ with $x \in \FU$ such that  $f|_{\FU}$ factors as $\FU \xto{g} \BA^{n}_{\FY} \xto{p} \FY$ where $g$ is an  \'etale  morphism and $p$ is the canonical projection. Applying \cite[Proposition 4.8]{AJP} we have that $f|_{\FU}$ and  $g$ are flat  morphisms and therefore, 
\[
\begin{array}{ccccc}
\dim_{x} f & = &\dim \widehat{\CO_{\FX,x} }\otimes_{\widehat{\CO_{\FY,y} }} k(y) &= &\dim \widehat{\CO_{\FX,x} } - \dim\widehat{\CO_{\FY,y} }\\
\dim_{x} g &=  &\dim \widehat{\CO_{\FX,x} }\otimes_{\widehat{\CO_{\BA^{n}_{\FY},g(x)}}} k(g(x))&= &\dim \widehat{\CO_{\FX,x} } - \dim\widehat{\CO_{\BA^{n}_{\FY},g(x)}}.
\end{array}
\]
Now, since $g$ is unramified  by   Corollary \ref{pnrdim0} we have that $\dim_{x} g =0$ and therefore $\dim_{x} f  = \dim\widehat{\CO_{\BA^{n}_{\FY},g(x)}}- \dim\widehat{\CO_{\FY,y} }=n$.
\end{proof}

\begin{propo} \label{ecppl}
Let $f:\FX \to \FY$ be a morphism of pseudo-finite type and let $\FX' \inc \FX$ be a  closed immersion given by the  ideal $\CI \subset \CO_{\FX}$ and put $f'=f|_{\FX'}$. If $f$ is smooth  at $x \in \FX'$, $n = \dim_{x} f$ and $y=f(x)$ the following conditions are equivalent:
\begin{enumerate}
\item
The morphism $f'$ is smooth  at $x$ and $\dim_{x} f'^{-1}(y)=n-m$.
\item
The natural sequence of $\CO_{\FX}$-modules
\[
0 \to \frac{\CI}{\CI^{2}} \to \om^{1}_{\FX/\FY} \otimes_{\CO_{\FX}}  \CO_{\FX'} \to \om^{1}_{\FX'/\FY}  \to 0
\]
is exact\footnote{Let $(X,\CO_{X})$ be a ringed space. We say that the sequence of  $\CO_{X}$-Modules \(0 \to \CF \to \CG \to \CH \to 0\) is exact at $x \in X$ if and only if $0 \to \CF_{x} \to \CG_{x} \to \CH_{x} \to 0$ is an exact sequence of $\CO_{X,x}$-modules.} at $x$ and, on a neighborhood of $x$, the displayed $\CO_{\FX'}$-Modules are locally free of  ranks $m,\, n$ and $n-m$, respectively.
\end{enumerate}
\end{propo}

\begin{proof}
Since $f:\FX \to \FY$ is a smooth morphism   at $x$, replacing $\FX$, if  necessary, by a smaller neighborhood  of  $x$, we may assume that $f: \FX=\spf(A) \to \FY=\spf(B)$ is a morphism in $\sfna$  smooth at $x$ and that $\FX'= \spf(A/I)$. Therefore, applying \cite[Proposition 4.8]{AJP} and  Corollary \ref{dimrango} we have that $\om^{1}_{\FX/\FY}$ is a locally free $\CO_{\FX}$-Module of  rank $n$. 

Let us prove that (1) $\Rightarrow $ (2). Replacing $\FX'$ with a smaller neighborhood of $x$ if necessary, we may also assume that $f':\FX' \to\FY$ is a smooth morphism. Then, by an argument along the lines of the previous paragraph, it follows that $\om^{1}_{\FX'/\FY}$ is a locally free $\CO_{\FX'}$-Module of  rank $n-m$. Zariski's Jacobian criterion for formal schemes (\cite[Corollary 4.15]{AJP}) implies that the sequence 
\[
0 \to \frac{\CI}{\CI^{2}} \to \om^{1}_{\FX/\FY} \otimes_{\CO_{\FX}}  \CO_{\FX'} \to \om^{1}_{\FX'/\FY}  \to 0
\]
is exact and split, from  which we deduce that $\CI/\CI^{2}$ is a locally free $\CO_{\FX'}$-Module of  rank $m$.

Conversely, applying \cite[(\textbf{0}, 5.5.4)]{EGA1} and 
\cite[Proposition 3.13]{AJP}
we deduce  that there exists an open formal subscheme $\FU \subset \FX'$ with $x \in \FU$ such that  
\[
0 \to \left(\frac{\CI}{\CI^{2}}\right)|_{\FU} \to (\om^{1}_{\FX/\FY} \otimes_{\CO_{\FX}}  \CO_{\FX'})|_{\FU} \to (\om^{1}_{\FX'/\FY})|_{\FU} \to 0
\]
is exact  and split. From  Zariski's Jacobian criterion it follows that $f'|_{\FU}$ is smooth  and therefore, $f'$ is smooth  at $x$.
\end{proof}

\begin{rem}
The natural sequence in Proposition \ref{ecppl} is  the Second Fundamental Exact  Sequence associated to the morphisms $\FX' \inc \FX \xto{f} \FX$ (\cite[Proposition 3.13]{AJP}).
\end{rem}

Locally, a pseudo-finite type morphism  $f: \FX \to \FY$  factors as  $\FU \overset{j} \inc \BD^{r}_{\BA^{s}_{\FY}} \xto{p} \FY$ where $j$ is a  closed immersion (see \cite[Proposition 1.7]{AJP}). In  Corollary \ref{criteriojacobiano} we provide a criterion in terms a matrix rank that tells whether  $\FU$ is smooth over $\FY$ or not.

\begin{parraf}
Let $\FY= \spf(A) \in\sfna$. Consider $\FX \subset \BD^{s}_{\BA^{r}_{\FY}}$ a closed formal subscheme given by an ideal $\CI=I^{\tr}$, with $I=\langle g_{1}\,, g_{2},\, \ldots,\, g_{k} \rangle \subset A\{\mathbf{T}\}[[\mathbf{Z}]]$ where $\mathbf{T}= T_{1},\, T_{2},\, \ldots,\, T_{r}$ and $\mathbf{Z}= Z_{1},\, Z_{2},\, \ldots,\, Z_{s}$ are two sets of of indeterminates. From \cite[3.14]{AJP} we have that  
\[\{\hd T_{1},\, \ldots ,\hd T_{r},\, \hd Z_{1},\,\ldots, \hd Z_{s}\}\]
is a basis of  $\om^{1}_{A\{\mathbf{T}\}[[\mathbf{Z}]]/A}$ and also that given $g \in A\{\mathbf{T}\}[[\mathbf{Z}]]$ it holds that:
\[
\hd g = \sum_{i=1}^{r}  \frac{\partial g}{\partial T_{i}} \hd T_{i} + \sum_{j=1}^{s}  \frac{\partial g}{\partial Z_{j}} \hd Z_{j},
\]
where $\hd$ is the complete canonical derivation of  $A\{\mathbf{T}\}[[\mathbf{Z}]]$ over $A$.
For any $g \in A\{\mathbf{T}\}[[\mathbf{Z}]]$, $w \in \{\hd T_{1},\,\ldots ,\hd T_{r},\, \hd Z_{1},\,\ldots, \hd Z_{s} \}$ and $x \in \FX$,  denote by $\frac{\partial g}{\partial w}(x) $ the image  of  $\frac{\partial g}{\partial w} \in A\{\mathbf{T}\}[[\mathbf{Z}]]$ in $k(x)$. We will call
\begin{equation*}
\Jac_{\FX/\FY}(x)=
    \begin{pmatrix}
	\frac{\partial g_{1}}{\partial T_{1}}(x) & \ldots & 
	\frac{\partial g_{1}}{\partial T_{r}}(x) & \frac{\partial g_{1}}{\partial Z_{1}}(x) & \ldots &	\frac{\partial g_{1}}{\partial Z_{s}}(x) \\
	\frac{\partial g_{2}}{\partial T_{1}}(x) & \ldots & 
	\frac{\partial g_{2}}{\partial T_{r}}(x) & \frac{\partial g_{2}}{\partial Z_{1}}(x) & \ldots &
	\frac{\partial g_{2}}{\partial Z_{s}}(x) \\
	\vdots  & \ddots & \vdots & \vdots  & \ddots  & \vdots\\
	\frac{\partial g_{k}}{\partial T_{1}}(x) & \ldots & 
	\frac{\partial g_{k}}{\partial T_{r}}(x) & \frac{\partial g_{k}}{\partial Z_{1}}(x) & \ldots &	\frac{\partial g_{k}}{\partial Z_{s}}(x) \\
	\end{pmatrix}
\end{equation*}
the \emph{Jacobian matrix of  $\FX$ over $ \FY$ at $x$}. This matrix depends on  the chosen generators of  $I$ and therefore, the notation  $\Jac_{\FX/\FY}(x)$ is not completely accurate. 
\end{parraf}

\begin{cor}{(Jacobian criterion for the affine formal space and the affine formal disc.)} \label{criteriojacobiano}
With the previous notations, the following assertions  are equivalent:
\begin{enumerate}
\item
The morphism $f: \FX  \to\FY$ is smooth  at $x$ and $\dim_{x} f= r+s-l$.
\item
There exists a subset $\{g_{1},\,, g_{2},\, \ldots,\, g_{l}\} \subset \{g_{1}\,, g_{2},\, \ldots,\, g_{k}\}$ such that  $\CI_{x}=\langle g_{1}\,, g_{2},\, \ldots,\, g_{l} \rangle\CO_{\FX,x}$ and $\rg(\Jac_{\FX/\FY}(x)) = l$.
\end{enumerate}
\end{cor}

\begin{proof}
Assume (1). By Proposition \ref{ecppl} we have that the sequence
\[
0 \to \frac{\CI}{\CI^{2}} \to \om^{1}_{ \BD^{s}_{\BA^{r}_{\FY}}/\FY} \otimes_{\CO_{ \BD^{s}_{\BA^{r}_{\FY}}}} \!\!\! \CO_{\FX} \to \om^{1}_{\FX/\FY}  \to 0
\]
is exact at $x$ and the corresponding $\CO_{\FX}$-Modules are locally free, in a neighborhood of $x$, of ranks $l$, $r+s$ and $r+s-l$, respectively. Therefore, 
\begin{equation} \label{sexcritjacobkev}
0 \to \frac{\CI}{\CI^{2}}\otimes_{\CO_{\FX}} k(x) \to   
\om^{1}_{ \BD^{s}_{\BA^{r}_{\FY}}/\FY} \otimes_{ \CO_{\BD^{s}_{\BA^{r}_{\FY}}}} \!\! k(x) \to  \om^{1}_{\FX/\FY}  \otimes_{\CO_{\FX}} k(x)  \to 0
\end{equation}
is an exact  sequence of  $k(x)$-vector spaces  of  dimension $l,\, r+s,\, r+s-l$, respectively. Thus, there exists a set $\{g_{1},\, g_{2},\, \ldots,\, g_{l}\} \subset \{g_{1},\, g_{2},\, \ldots,\, g_{k}\}$ such that  $\{g_{1}(x),\, g_{2}(x),\, \ldots,\, g_{l}(x)\}$ provides a basis of  $\CI/\CI^{2}\otimes_{\CO_{\FX}} k(x)$  at $x$. By Nakayama's lemma it holds that $\CI_{x}=\langle g_{1}\,, g_{2},\, \ldots,\, g_{l} \rangle\CO_{\FX,x}$. Besides, from  the exactness of  the sequence (\ref{sexcritjacobkev}) and from  the equivalence of categories \cite[(10.10.2)]{EGA1} we deduce that the set 
\[\{\hd g_{1}(x),\, \hd g_{2}(x),\, \ldots,\, \hd g_{l}(x)\} \subset \om^{1}_{ A\{\mathbf{T}\}[[\mathbf{Z}]]/A} \otimes_{A\{\mathbf{T}\} [[\mathbf{Z}]]} k(x)\] 
is linearly independent. Therefore, $\rg(\Jac_{\FX/\FY}(x)) = l$. 

Conversely, from  the Second Fundamental Exact Sequence associated to the morphisms $\FX \inc \BD^{s}_{\BA^{r}_{\FY}} \to \FY$ \cite[Proposition 3.13]{AJP} we get the exact sequence
\[
 \frac{\CI}{\CI^{2}}\otimes_{\CO_{\FX}} k(x) \to   
\om^{1}_{ \BD^{s}_{\BA^{r}_{\FY}}/\FY} \otimes_{ \CO_{\BD^{s}_{\BA^{r}_{\FY}}}} \!\! k(x) \to  \om^{1}_{\FX/\FY}  \otimes_{\CO_{\FX}} k(x)  \to 0.
\]
Since  $\rg(\Jac_{\FX/\FY}(x)) = l$,   we have  that 
\[
\{\hd g_{1}(x),\,, \hd g_{2}(x),\, \ldots,\, \hd g_{l}(x)\} \subset \om^{1}_{ A\{\mathbf{T}\}[[\mathbf{Z}]]/A} \otimes_{A\{\mathbf{T}\} [[\mathbf{Z}]]} k(x)
\]
is a linearly independent set. Extending this set to a basis of the vector space $\om^{1}_{ A\{\mathbf{T}\}[[\mathbf{Z}]]/A} \otimes_{A\{\mathbf{T}\} [[\mathbf{Z}]]} k(x)$, by Nakayama's lemma we find a basis 
\(\mathcal{B} \subset \om^{1}_{ A\{\mathbf{T}\}[[\mathbf{Z}]]/A}\)
 such that \(\{\hd g_{1},\, \hd g_{2},\, \ldots, \, \hd g_{l}\} \subset \mathcal{B}\)
and therefore  
\[\{\hd g_{1},\,, \hd g_{2},\, \ldots,\, \hd g_{l}\} \subset\om^{1}_{ A\{\mathbf{T}\}[[\mathbf{Z}]]/A} \otimes_{A\{\mathbf{T}\} [[\mathbf{Z}]]} A\{\mathbf{T}\}[[\mathbf{Z}]]/I\]
is a linearly independent set  at $x$. Thus the  set $\{g_{1}\,, g_{2,},\, \ldots,\, g_{l}\}$ provides a basis  of  $\CI/\CI^{2}$ at $x$ and by the equivalence of categories \cite[(10.10.2)]{EGA1} we have that the sequence of $\CO_\FX$-Modules
\[
0 \to \frac{\CI}{\CI^{2}} \to \om^{1}_{ \BD^{s}_{\BA^{r}_{\FY}}/\FY} \otimes_{\CO_{ \BD^{s}_{\BA^{r}_{\FY}}}} \!\!\! \CO_{\FX} \to \om^{1}_{\FX/\FY}  \to 0
\]
is split exact at $x$ of locally free Modules  of  ranks $l$, $r+s$ and $r+s-l$, respectively. Applying Proposition \ref{ecppl} it follows that $f$ is smooth  at $x$ and $\dim_{x} f= r+s-l$.
\end{proof}

Notice that the matrix form of the Jacobian criterion for the  affine formal space  and the affine formal disc (Corollary \ref{criteriojacobiano}) generalize the usual matrix form of the Jacobian criterion for the  affine space in $\sch$ (\cite[Ch. VII, Theorem (5.14)]{at}).


\section{\'Etale morphisms} \label{sec5}
The  main results of this section are consequences of those obtained in Sections \ref{sec3} and \ref{sec4}. They will allow us to characterize in Section \ref{sec6} two important classes  of \'etale morphisms: open immersions and completion morphisms.

\begin{propo} \label{efn}
Given $f:\FX \to \FY$ in $\sfn$  let $\CJ \subset \CO_{\FX}$ and $\CK \subset \CO_{\FY}$ be $f$-compatible ideals of  definition. Using these ideals, set
\[f=\dirlim {n \in \NN} f_{n}.\] If $f_{n}:X_{n} \to Y_{n}$ is \'etale, $\forall n \in \NN$, then $f$ is \'etale.
\end{propo}

\begin{proof}
The sum of Proposition \ref{nrfn} and Proposition \ref{lfn}.
\end{proof}

\begin{cor} \label{eadfn}
Let $f:\FX \to \FY$ be an  \emph{adic} morphism  in $\sfn$ and let $\CK \subset \CO_{\FY}$ be an ideal of  definition. Consider $\{f_{n}\}_{n \in \NN}$ the  direct system  of  morphisms of  schemes associated to the ideals of  definition $\CK \subset \CO_{\FY}$ and $\CJ=f^{*}(\CK)\CO_{\FX} \subset \CO_{\FX}$.
The morphism $f$ is \'etale if and only if the morphisms $f_{n}:X_{n} \to Y_{n}$ are \'etale $\forall n \in \NN$.
\end{cor}
\begin{proof}
It follows from Proposition  \ref{nrfn} and  Corollary \ref{ladfn}.
\end{proof}

\begin{propo} \label{fetf0et}
Let $f:\FX \to \FY$ be an  \emph{adic} morphism  in $\sfn$ and let $f_{0}:X_{0} \to Y_{0}$ be the morphism of schemes associated to the ideals of definition $\CK \subset \CO_{\FY}$ and $\CJ =f^{*}(\CK)\CO_{\FX} \subset \CO_{\FX}$. Then, $f$ is \'etale if and only if $f$ is flat and $f_{0}$ is \'etale.
\end{propo}

\begin{proof}
Put together Proposition \ref{fnrf0nr} and Corollary  \ref{flf0l}.
\end{proof}

Note that Example \ref{peynopen}  shows that in the non adic case the last two results do not hold  and also that, in general, the  converse of  Proposition \ref{efn} is not true.  

\begin{propo}
Let $f$ be a pseudo-finite type morphism in $\sfn$ and choose $\CJ\subset \CO_{\FX}$ and $\CK \subset \CO_{\FY}$ $f$-compatible ideals of definition. Write 
\[f=\dirlim {n \in \NN} f_{n}.\] 
If $f_{0}:X_{0} \to Y_{0}$ is \'etale, $f$ is flat and  $f^{-1}(f(x))= f_{0}^{-1}(f(x))$, for all $x \in \FX$, then $f$ is \'etale.
\end{propo}

\begin{proof}
It follows from Corollary  \ref{corf0imfpnr} and Corollary \ref{corf0imfpl}.
\end{proof}

Example \ref{peynopen} shows that the  converse of the last result  is not true. Next Proposition gives us a local characterization of  \'etale morphisms.

\begin{propo} \label{caractlocalpe}
Let $f:\FX \to \FY$ be a morphism in $\sfn$ of  pseudo-finite type, let  $x \in \FX$ and $y = f(x)$, the following conditions are equivalent:
\begin{enumerate}
\item [(1)]
 $f$ is \'etale  at $x$.
\item [(2)]
$\CO_{\FX,x}$ is a  formally \'etale $\CO_{\FY,y}$-algebra for the adic topologies.
\item [$(2')$]
$\widehat{\CO_{\FX,x}}$ is a formally \'etale $\widehat{\CO_{\FY,y}}$-algebra  for the adic topologies.
\item [(3)]
$f$ is flat at $x$ and $f^{-1}(y)$ is a $k(y)$-formal scheme \'etale  at $x$.
\item [(4)]
$f$ is flat and unramified  at $x$.
\item [$(4')$]
$f$ is flat at $x$ and $(\om^{1}_{\FX/\FY})_{x}=0$. 
\item [(5)]
$f$ is smooth  at $x$ and a quasi-covering at $x$.
\end{enumerate}
\end{propo}

\begin{proof}
Applying Proposition \ref{caraclocalpnr} and Proposition \ref{pligplafibr} we have  that
\[(5)  \Leftarrow (1) \dimp (2) \dimp (2') \dimp (3) \imp (4) \dimp (4').\]
Let $C := \widehat{\CO_{\FX,x}} \otimes_{\widehat{\CO_{\FY,y}}} k(y)$.
To show  (4) $\imp$ (5),  by  Corollary \ref{corpnrimplcr} it is only left to prove that $f$ is smooth  at $x$. By hypothesis, we have that $f$ is unramified  at $x$ and by Proposition \ref{caraclocalpnr},  it follows that $C = k(x)$ and $k(x)|k(y)$ is a  finite separable extension, therefore, formally \'etale. Since $f$ is flat at $x$, by Proposition \ref{pligplafibr} we conclude that $f$ is smooth  at $x$. 

To prove that (5) $\imp$ (1), it suffices to check that $f$ is unramified  at $x$ or, equivalently by Proposition \ref{caraclocalpnr}, that  $C = k(x)$ and that $k(x)|k(y)$ is a finite separable extension. As $f$ is smooth at $x$, applying Proposition \ref{pligplafibr}, we have that $\widehat{\CO_{\FX,x}}$ is a formally smooth $\widehat{\CO_{\FX,x}}$-algebra  for the adic topologies. Then by  base-change it holds that $C$ is a  formally smooth $k(y)$-algebra. By \cite[(\textbf{0}, 19.3.8)]{EGA41} we have that $C$ is a formally smooth $k(y)$-algebra  for the topologies given by the maximal ideals and from \cite[Lemma 1, p. 216]{ma2} it holds that $C$ is a regular local ring. Besides, by hypothesis we have  that  $C$ is a finite $k(y)$-module, therefore, an artinian ring, so $C = k(x)$. Since   $k(x) = C$ is a formally smooth $k(y)$-algebra  we have that $k(x)|k(y)$ is a separable extension (\emph{cf.} \cite[(\textbf{0}, 19.6.1)]{EGA41}).  
\end{proof} 

\begin{cor} \label{corcaractlocalpe}
Let $f:\FX \to \FY$ be a  pseudo-finite type morphism in $\sfn$. The following conditions are equivalent:
\begin{enumerate}
\item [(1)]
 $f$ is \'etale.
\item [(2)]
For all $x \in \FX$, $\CO_{\FX,x}$ is a formally \'etale $\CO_{\FY,f(x)}$-algebra for  the adic topologies.
\item [$(2')$]
For all $x \in \FX$, $\widehat{\CO_{\FX,x}}$ is a formally \'etale  $\widehat{\CO_{\FY,f(x)}}$-algebra for  the adic topologies.
\item  [(3)]
For all $x \in \FX$, $f^{-1}(f(x))$ is a $k(f(x))$-formal scheme \'etale  at $x$ and $f$ is flat. 
 \item [(4)]
$f$ is flat and unramified.
\item [$(4')$]
$f$ is flat  and $\om^{1}_{\FX/\FY}=0$.
\item [(5)]
$f$ is smooth  and a quasi-covering.

\end{enumerate}

\end{cor}

\begin{ex}\label {pcf+plnope}
Given a field $K$, the  canonical morphism  $\BD^{1}_{K} \to \spec(K)$ is smooth,  pseudo-quasi-finite  but  it is not  \'etale.
\end{ex}

In $\sch$ a morphism is \'etale if and only if it is smooth and quasi-finite. The  previous example shows that in $\sfn$ there are smooth  and pseudo-quasi-finite morphisms that are not \'etale. That is why we consider quasi-coverings in $\sfn$ (see Definition \ref{defncuasireves}) as the right generalization of quasi-finite morphisms in $\sch$.

\section[Structure theorems]{Structure theorems of the infinitesimal lifting properties} \label{sec6}

We begin with two results that will be used in the proof  of  the remainder results of  this  section.
  
\begin{propo} \label{levantpeusform}
Consider a formally \'etale morphism  $f:\FX \to \FY$ and a morphism $g\colon \FS \to \FY$, both in $\sfn$. Take $\CL \subset \CO_{\FS}$ an ideal of definition of $\FS$ and write
\[\FS = \dirlim {n \in \NN} S_{n}.\] 
If  $h_{0}: S_{0} \to \FX$ is a morphism in $\sfn$ that makes the  diagram 
\begin{diagram}[height=2em,w=2em,p=0.3em,labelstyle=\scriptstyle]
S_{0}         & \rTinc    & \FS	\\
\dTto^{h_{0}} &           & \dTto_{g}	\\ 
\FX           & \rTto^{f} & \FY\\  
\end{diagram}
commutative, where $S_{0} \inc \FS$ is the  canonical  closed immersion, then
there exists a unique $\FY$-morphism $l: \FS \to \FX$  in $\sfn$ such that  $l|_{S_{0}}=h_{0}$.
\end{propo}

\begin{proof}
By induction on $n$ we are going to construct a collection of  morphisms $\{h_{n}: S_{n} \to \FX\}_{n \in \NN}$ such that the diagrams
\begin{diagram}[height=2.2em,w=2.5em,labelstyle=\scriptstyle]
S_{n-1}	&        &       &       &\\
		&\rdTto	\rdTto(4,2) \rdTto(2,4)_{h_{n-1}}	&    &      &   \\
        &       &S_{n}	&\rTinc & \FS\\
        &       & \dTto_{h_{n}} &	& \dTto_{g}\\
        &       & \FX	&\rTto^{f}&	\FY\\
\end{diagram}
commute.
For $n =1$, by \cite[2.4]{AJP} there exists a unique morphism $h_{1}: S_{1}  \to \FX$ such that  $h_{1} |_{S_{0}}= h_{0}$ and $g|_{S_{1}}= f \circ h_{1}$. Now let $n \in \NN, \,n > 1$ and suppose we already have for all $0< k <n$  morphisms $h_{k}: S_{k} \to \FX$ such that $h_{k} |_{S_{k-1}}= h_{k-1}$ and $g|_{S_{k}}= f \circ h_{k}$. Then by \cite[\emph{loc. cit.}]{AJP} there exists a unique morphism $h_{n}: S_{n} \to \FX$ such that $h_{n} |_{S_{n-1}}= h_{n-1}$ and $g|_{S_{n}}= f \circ h_{n}$. It is straightforward that \[l:= \dirlim {n \in \NN} h_{n}\] is a morphism of  formal schemes and is the  unique one such that the  diagram
\begin{diagram}[height=2em,w=2em,p=0.3em,labelstyle=\scriptstyle]
S_{0}         & \rTinc     & \FS	\\
\dTto^{h_{0}} & \ldTto_{h} & \dTto_{g}	\\ 
\FX           & \rTto_f   & \FY \\  
\end{diagram}
commutes.
\end{proof}


\begin{cor} \label{iso0peis}
Let $f:\FX \to \FY$ be an   \'etale morphism in $\sfn$ and  $\CJ \subset \CO_{\FX}$ and $\CK \subset \CO_{\FY}$ $f$-compatible ideals of  definition   such that the corresponding morphism $f_{0}:X_{0} \to Y_{0}$ is an isomorphism. Then $f$ is an isomorphism. 
\end{cor}

\begin{proof}
By Proposition \ref{levantpeusform} there exists a (unique) morphism $g: \FY \to \FX$ such that  the  following diagram is commutative
\begin{diagram}[height=2em,w=2.5em,p=0.3em,labelstyle=\scriptstyle]
Y_{0}              & \rTinc          & \FY \\
\dTto^{f_{0}^{-1}} & \ldTto(2,4)^{g} &     \\
X_{0}              &                 & \dTto_{1_{\FY}}\\
\dTinc             &                 & \\
\FX                & \rTto^{f}       & \FY\\  
\end{diagram}
Then,  by  \cite[Proposition 2.13]{AJP} it follows that $g$ is an \'etale morphism. Thus, applying  Proposition \ref{levantpeusform} we have that there exists a (unique) morphism $f': \FX \to \FY$ such that   the  following diagram is commutative
\begin{diagram}[height=2em,w=2.5em,p=0.3em,labelstyle=\scriptstyle]
X_{0}	      &   \rTinc	            & \FX \\
\dTto^{f_{0}} &   \ldTto(2,4)^{f'}   &	\\
Y_{0}         &             & \dTto_{1_{\FX}}	\\
 \dTinc       &             & \\
\FY  	      & \rTto^{g}   &\FX\\  
\end{diagram}
From $f \circ g=1_{\FY}$ and $g \circ f'=1_{\FX}$ we deduce that $f=f'$ and therefore $f$  is an isomorphism.
\end{proof}

In $\sch$ open immersions are characterized as being those \'etale morphisms that are radicial (see \cite[(17.9.1)]{EGA44}). In the  following theorem we  extend this characterization and relate open immersions in formal schemes  with their counterparts in schemes.
 
\begin{thm} \label{caractencab}  
Let $f:\FX \to \FY$ be a morphism in $\sfn$. The following conditions are equivalent:
\begin{enumerate}
\item
$f$ is an open immersion.
\item
$f$ is adic, flat and if $\CK \subset \CO_{\FY}$ is an ideal of  definition such that $\CJ= f^{*}(\CK) \CO_{\FX} \subset \CO_{\FX}$,  the  associated  morphism of  schemes $f_{0}:X_{0} \to Y_{0}$ is an open immersion.
\item
$f$ is \emph{adic} \'etale and radicial.
\item
There are $\CJ \subset \CO_{\FX}$ and $\CK \subset \CO_{\FY}$ $f$-compatible ideals of  definition such that the morphisms $f_{n}:X_{n} \to Y_{n}$ are open immersions, for all $n \in \NN$.
\end{enumerate}
\end{thm}

\begin{proof}
The  implication  (1) $\Rightarrow$ (2) is immediate. Given $\CK \subset \CO_{\FY}$ an ideal of  definition, assume (2)  and let us show (3). Since $f_{0}$ is an open immersion, is radicial, so,  $f$ is radicial (see Definition \ref{rad} and its attached paragraph). Furthermore, $f$ is flat and $f_{0}$ is an \'etale morphism then $f$ is \'etale (see Proposition \ref{fetf0et}). 
Let us prove that (3) $\Rightarrow$ (4). Given $\CK \subset \CO_{\FY}$ an ideal of  definition and $\CJ = f^{*}(\CK) \CO_{\FX}$, by  Corollary \ref{eadfn} the morphisms $f_{n}:X_{n} \to Y_{n}$ are  \'etale, for all $n \in \NN$. The morphisms $f_{n}$ are also radicial for all $n \in \NN$ (see Definition \ref{rad}) and thus by \cite[(17.9.1)]{EGA44} it follows that $f_{n}$ is an open immersion, for each $n \in \NN$. Finally, suppose  that  (4) holds and let us see that $f$ is an open immersion. With the notations in (4), there exists an open subset $U_{0} \subset Y_{0}$ such that  $f_{0}$ factors as
\[
X_{0} \xto{f'_{0}} U_{0} \overset{i_{0}} \inc Y_{0}
\]
where $f'_{0}$ is an isomorphism and $i_{0}$ is the canonical inclusion. Let $\FU \subset \FY$ be the  open formal subscheme with underlying topological space $U_{0}$. Since the  open immersion $i: \FU \to \FY$ is \'etale, then Proposition \ref{levantpeusform} implies that there exists a morphism $f': \FX \to \FU$ of formal schemes such that  the  diagram
\begin{diagram}[height=1.5em,w=2em,labelstyle=\scriptstyle]
\FX	    &       & \rTto^{\qquad f}&			&	   &\FY\\
\uTinc	& \rdTto(3,2)^{f'}&	&  		&\ruTinc^{i}&	\uTinc \\
		&	    &            & \FU&  & \\
		&	    &            &\uTto	&	   & 	          \\
X_{0}	&       & \hLine^{\qquad f_{0}}&\VonH	&\rTto   &Y_{0}\\
		& \rdTto(3,2)_{f'_{0}}       &	&	  &\ruTinc^{i_{0}} &    \\
		&	    &		& U_{0}  				&	 &   \\
\end{diagram} 
is commutative.
Since the morphisms $f_{n}$ are \'etale, for all $n \in \NN$, Proposition \ref{efn} implies that $f$ is \'etale. By \cite[Proposition 2.13]{AJP} we have that $f'$ is \'etale  and applying  Corollary \ref{iso0peis},  $f'$ is an isomorphism and therefore, $f$ is an open immersion.
\end{proof}

\begin{cor}\label{diagopen}
Let $f\colon\FX \to \FY$ be a pseudo-finite type morphism in $\sfn$. Then $f$ is unramified if and only if the diagonal morphism $\Delta_f \colon \FX \to \FX \times_\FY \FX$ is an open embedding.
\end{cor}

\begin{proof}
Take $\CJ \subset \CO_{\FX}$ and $\CK \subset \CO_{\FY}$ $f$-compatible ideals of definition such that $f$ can be expressed as the limit of maps of usual schemes $f_n \colon X_n \to Y_n$, $n \in \NN$.
The morphism $f\colon\FX \to \FY$ is unramified if and only if $f_n$ is unramified for all $n \in \NN$ by Proposition \ref{nrfn}. By \cite[Corollaire (17.4.2)]{EGA44} this is equivalent to $\Delta_{f_n}: X_n \to X_n \times_{Y_n} X_n$ being an open embedding for all $n \in \NN$. But this, in turn, is equivalent to the fact that $\Delta_f \colon \FX \to \FX \times_\FY \FX$ is an open embedding by Theorem \ref{caractencab}.
\end{proof}

Every completion morphism is a pseudo-closed immersion that is flat (\emph{cf.} Proposition \ref{caractcom}). Next, we prove that this condition is also sufficient. Thus, we obtain a criterion to determine whether a $\FY$-formal scheme $\FX$ is the completion of  $\FY$ along  a closed formal subscheme.
 
\begin{thm} \label{caracmorfcompl}
Let $f:\FX \to \FY$ be a morphism in $\sfn$ and let $\CJ \subset \CO_{\FX}$ and $\CK \subset \CO_{\FY}$ be $f$-compatible ideals of  definition. Let $f_{0}:X_{0} \to Y_{0}$ be the corresponding morphism of ordinary schemes. The following conditions are equivalent:
\begin{enumerate}
\item
There exists a closed formal subscheme $\FY' \subset \FY$ such that $\FX = \FY_{/\FY'}$ and $f$ is the morphism of completion of $\FY$ along $\FY'$.
\item
The morphism $f$ is a flat pseudo-closed immersion.
\item
The morphism $f$ is  \'etale  and $f_{0}:X_{0} \to Y_{0}$ is a  closed immersion.
\item
The morphism $f$ is a smooth pseudo-closed immersion.
\end{enumerate}
\end{thm}

\begin{proof}
The implication (1) $\imp$ (2) is Proposition   \ref{caractcom}. Let us show that (2) $\imp$ (3). Since $f$ is a pseudo-closed immersion, by  Corollary \ref{pecigf0ecnr}  we have that $f$ is unramified. Then as  $f$ is flat,  Corollary \ref{corcaractlocalpe} establishes that $f$ is \'etale. The equivalence (3) $\dimp$ (4) is consequence of Corollary \ref{pecigf0ecnr}.  Finally, we show that (3) $\imp$ (1). By hypothesis, the   morphism $f_{0}: X_{0} \to \FY$ is a  closed immersion. Consider  $\kappa: \FY_{/ X_{0}} \to \FY$ the  morphism of  completion of  $\FY$ along  $X_{0}$ and let us prove that $\FX $ and $\FY_{/ X_{0}}$ are $\FY$-isomorphic. By Proposition  \ref{caractcom} the  morphism $\kappa$ is \'etale  so,  applying Proposition \ref{levantpeusform}, we have that there exists a $\FY$-morphism $\varphi: \FX \to \FY_{/ X_{0}}$ such that  the  following diagram is commutative
\begin{diagram}[height=1.5em,w=2em,labelstyle=\scriptstyle]
\FX    &                      &  \rTto^{\qquad f} & &  &\FY\\
\uTinc & \rdTto(3,2)^{\varphi}&	&  		             &\ruTto^{\kappa}& \uTinc \\
       &                      && \FY_{/ X_{0}}      &  & \\
       &                      && \uTinc             &  & \\
X_{0}  & &\hLine^{\qquad f_{0}}& \VonH       &\rTto    & Y_{0} \\
       & \rdTto(3,2)_{\varphi_{0}=1_{X_{0}}} &&     &\ruTinc^{f_{0}} & \\
       &                      && X_{0}              & &   \\
\end{diagram}
From \cite[Proposition 2.13]{AJP} it follows that $\varphi$ is \'etale  and then by  Corollary \ref{iso0peis} we get that $\varphi$ is an isomorphism.
\end{proof}

\begin{rem}
A consequence of the proof of (3) $\imp$ (1) is the following: Given $\FY$ in $\sfn$ and a closed formal subscheme $\FY' \subset \FY$ defined by the  ideal $\CI \subset \CO_{\FY}$,  then for every ideal of  definition $\CK \subset \CO_{\FY}$  of  $\FY$, it holds that 
\[\FY_{/\FY'} = \FY_{/Y'_{0}}\] 
where $Y'_{0} =  (\FY',\CO_{\FY}/(\CI+\CK))$.
\end{rem}


\begin{parraf}
Given  a scheme $Y$ and a closed subscheme $Y_{0} \subset Y$ with the  same topological space, the  functor $X \leadsto X \times_{Y} Y_{0}$ defines an equivalence between the  category of   \'etale $Y$-schemes and the category of   \'etale $Y_{0}$-schemes by \cite[(18.1.2)]{EGA44}. In the  next theorem we extend this equivalence to the category of  locally noetherian formal schemes.
A special case of this theorem, namely when $\FY$ is smooth over a noetherian ordinary base scheme, appears in \cite[Proposition 2.4]{y}.
\end{parraf}

\begin{propo} \label{teorequivet}
Let $\FY$ be in $\sfn$ and $\CK \subset \CO_{\FY}$ an ideal of  definition such that 
\[\FY = \dirlim {n \in \NN} Y_{n}.\] Then the  functor
\[
\begin{matrix}
\textrm { \'etale adic }\FY\textrm{-formal schemes}& \xto{F} &\textrm{\'etale }Y_{0}\textrm{-schemes}\\
\FX& \leadsto& \FX \times_{\FY} Y_{0}\\
\end{matrix}
\]
is an equivalence of  categories.
\end{propo}

\begin{proof}
By \cite[IV, \S 4, Theorem 1]{mcl} it suffices to prove that: (a) $F$ is full and faithful; and (b) Given $X_{0}$ an  \'etale $Y_{0}$-scheme  there exists an  \'etale adic $\FY$-formal scheme $\FX$ such that  $F(\FX)=\FX \times_{\FY}Y_{0} \cong X_{0}$.

The assertion (a) is an immediate consequence of  Proposition \ref{levantpeusform}.

Let us show (b). Given  $X_{0}$ an \'etale $Y_{0}$-scheme  in $\sch$ by \cite[(18.1.2)]{EGA44} there exists $X_{1}$ a locally noetherian \'etale $Y_{1}$-scheme   such that   $X_{1}\times_{Y_{1}}Y_{0}  \cong X_{0}$. Reasoning by induction on $n \in \NN$ and using \cite[\emph{loc. cit.}]{EGA44}, we get a family of schemes $\{X_{n}\}_{n \in \NN}$ such that, for each $n \in \NN$, $X_{n}$ is a locally noetherian \'etale $Y_{n}$-scheme and  $X_{n}\times_{Y_{n}}Y_{n-1}  \cong X_{n-1}$, for $n > 0$. Then 
\[\FX:= \dirlim {n \in \NN} X_{n}\]
is a locally noetherian adic $\FY$-formal scheme (by \cite[(10.12.3.1)]{EGA1}), 
\[\FX \times_{\FY}Y_{0} \underset{\textrm{\cite[(10.7.4)]{EGA1}} }= \dirlim {n \in \NN} (X_{n} \times_{Y_{n}} Y_{0}) =X_{0}\] and 
$\FX$ is an  \'etale $\FY$-formal scheme (see Proposition \ref{efn}).
\end{proof}

\begin{rem}
It seems plausible that there is a theory of an algebraic fundamental group for formal schemes that classifies \emph{adic} \'etale surjective maps onto a noetherian formal scheme $\FX$. If this is the case, the previous theorem would imply that it agrees with the fundamental group of $X_0$. We also consider feasible the existence of a bigger fundamental group classifying arbitrary \'etale surjective maps onto a noetherian formal scheme $\FX$, that would give additional information on $\FX$.
\end{rem}

\begin{cor}
Let $f:\FX \to \FY$ be an \'etale morphism  in $\sfn$. Given $\CJ \subset \CO_{\FX},$ and $\CK \subset \CO_{\FY}$ $f$-compatible ideals of  definition, if the induced  morphism $f_{0}:X_{0} \to Y_{0}$ is \'etale, then $f$ is adic \'etale.
\end{cor}

\begin{proof}
By Proposition \ref{teorequivet} there is an adic \'etale morphism $f':\FX' \to\FY$ in $\sfn$  such that  $\FX' \times_{\FY} Y_{0}= X_{0}$. Therefore by Proposition \ref{levantpeusform} there exists a morphism of formal schemes $g: \FX \to \FX'$ such that  the  diagram
\begin{diagram}[height=1.5em,w=2em,labelstyle=\scriptstyle]
\FX	    &      &          & \rTto^{f} &  &\FY\\
\uTinc  & \rdTto(3,2)^{g} &&          & \ruTto^{f'} & \uTinc \\
		&      &           & \FX'	   &  & \\
		&      &           & \uTto	   &  & \\
X_{0}	&      & \hLine    & \VonH	   &\rTto^{f_{0}} & Y_{0} \\
		& \rdTto(3,2)_{g_{0}=1_{X_{0}}} & &	 & \ruTto^{f'_{0}} & \\
		&      &           & X_{0} &  &  \\
\end{diagram} 
is commutative.
Applying \cite[Proposition 2.13]{AJP} we have that $g$ is \'etale  and from Corollary \ref{iso0peis} we deduce that $g$ is an isomorphism and therefore, $f$ is adic \'etale.
\end{proof}

\begin{cor} \label{fnetfet}
Let $f:\FX \to \FY$ be a morphism in $\sfn$. The  morphism $f$ is adic \'etale if and only if there exist $f$-compatible ideals of definition $\CJ \subset \CO_{\FX}$ and $\CK \subset \CO_{\FY}$ such that the induced morphisms   $f_{n}:X_{n} \to Y_{n}$ are \'etale, for all $n \in \NN$.
\end{cor}

\begin{proof}
If $f$ is adic \'etale, given $\CK \subset \CO_{\FY}$ an ideal of  definition, take $\CJ=f^{*}(\CK) \CO_{\FX}$ the corresponding ideal of definition of  $\FX$. By  base change, we have that the morphisms  $f_{n}:X_{n} \to Y_{n}$ are  \'etale, for all $n \in \NN$. The  converse is a consequence of Proposition \ref{efn} and of the previous Corollary.
\end{proof}

Proposition \ref{teorequivet} says that given \[\FY= \dirlim {n \in \NN} Y_{n}\] in $\sfn$ and $X_{0}$ an  \'etale $Y_{0}$-scheme there exists  a unique (up to isomorphism)  \'etale $\FY$-formal scheme $\FX$ such that  $\FX \times_{\FY} Y_{0}=X_{0}$.  But, what happens when $X_{0}$ is a smooth $Y_{0}$-scheme? 

\begin{propo} \label{equivlocal}
Let $\FY$ be in $\sfn$ and with respect to an ideal of  definition $\CK \subset \CO_{\FY}$ let us write 
\[\FY = \dirlim {n \in \NN} Y_{n}.\] 
Given $f_{0}:X_{0} \to Y_{0}$ a  morphism in $\sch$  smooth at $x \in X_{0}$, there exists an open subset  $U_{0} \subset X_{0}$, with $x \in U_{0} $ and a smooth adic  $\FY$-formal scheme   $\FU$  such that  $\FU \times_{\FY} Y_{0} \cong U_{0}$.
\end{propo}

\begin{proof}
 Since this is a local question in $\FY$, we may assume that $\FY= \spf(B)$ is in $\sfna$, $\CK = K^{\tr}$ with $K \subset B$ an  ideal of  definition of the adic ring $B$, $B_{0} = B/K$ and that $f_{0}:X_{0}= \spec(A_{0}) \to Y_{0} = \spec(B_{0})$ is a morphism in $\scha$ smooth at $x \in X_{0}$.  By 
 Proposition \ref{factpl} there exists an open subset $U_{0} \subset X_{0}$ with $x \in U_{0}$ such that  $f_{0} |_{U_{0}}$ factors as 
 \[
 U_{0} \xto{f'_{0}} \BA_{Y_{0}}^{n}= \spec(B_{0}[\mathbf{T}]) \xto{p_{0}} Y_{0}
 \]
where $\mathbf{T}= T_{1},\, T_{2},\, \ldots,\, T_{r}$ is a set of indeterminates,  $f'_{0}$ is an \'etale morphism and $p_{0}$ is the canonical projection. The  morphism $p_{0}$ lifts to a projection morphism $p: \BA_{\FY}^{n}= \spf(B\{\mathbf{T}\}) \to \FY$ such that the square in the following diagram is cartesian
 \begin{diagram}[height=2em,w=2em,labelstyle=\scriptstyle]
      &                & \BA_{\FY}^{n}   & \rTto^{p}     & \FY\\
      &                & \uTinc          &               & \uTinc \\
U_{0} & \rTto^{f'_{0}} & \BA_{Y_{0}}^{n} & \rTto^{p_{0}} & Y_{0}\\
\end{diagram}
Applying Proposition \ref{teorequivet}, there exists a locally noetherian \'etale adic $\BA_{\FY}^{n}$-formal scheme $\FU$ such that  $U_{0} \cong \FU \times_{\BA_{\FY}^{n}} \BA_{Y_{0}}^{n}$. Then $\FU$ is an smooth adic $\FY$-formal scheme  such that  $U_{0} \cong \FU \times_{\FY} Y_{0}$.
\end{proof}


The next theorem transfers the local description of unramified morphisms known in the case of schemes (\cite[(18.4.7)]{EGA44}) to the framework  of  formal schemes.

\begin{thm} \label{tppalnr}
Let $f:\FX \to \FY$ be a morphism in $\sfn$ unramified  at $x \in \FX$. Then there exists an open subset $\FU \subset \FX$ with $x \in \FU$ such that  $f|_{\FU}$ factors as
\[
\FU \xto{\kappa} \FX' \xto{f'} \FY
\]
where $\kappa$ is a pseudo-closed immersion  and $f'$ is an \emph{adic} \'etale morphism.
\end{thm}

\begin{proof}
Let $\CJ \subset \CO_{\FX}$ and $\CK \subset \CO_{\FY}$ be ideals of definition.
The morphism of  schemes $f_{0}$ associated to these ideals is unramified at $x$ (Proposition \ref{nrfn}) and by \cite[(18.4.7)]{EGA44} there exists an open set  $U_{0}\subset X_{0}$ with $x \in U_{0}$ such that  $f_{0}|_{U_{0}}$ factors as
\[
U_{0}\rTinc^{\kappa_{0}} X'_{0} \rTto^{f'_{0}} Y_{0}
\]
where $\kappa_{0}$ is a  closed immersion and $f'_{0}$ is an \'etale morphism. Proposition \ref{teorequivet} implies that there exists  an  \'etale adic morphism $f':\FX' \to\FY$ in $\sfn$ such that  $\FX' \times_{\FY} Y_{0} = X'_{0}$. Now, if $\FU \subset \FX$ is the  open formal scheme  with underlying topological space  $U_{0}$, by  Proposition \ref{levantpeusform} there exists a morphism $\kappa: \FU \to \FX'$ such that  the  following diagram commutes
\begin{diagram}[height=1.5em,w=2em,labelstyle=\scriptstyle]
\FU    &  & \rTto^{f|_{\FU}}  &	       &           & \FY \\
\uTinc & \rdTto(3,2)^{\kappa} &&      &\ruTto^{f'}& \uTinc \\
       &  &                   & \FX'  &           & \\
       &  &                   &\uTto  &           & \\
U_{0}  &  &\hLine^{f_{0}|_{U_{0}}}&\VonH	&\rTto & Y_{0} \\
       & \rdTinc(3,2)_{\kappa_{0}}  &&& \ruTto^{f'_{0}} & \\
       &  &                   & X'_{0} &	        & \\
\end{diagram}
Since $f$ is unramified,  by \cite[Proposition 2.13]{AJP} it holds that $\kappa$ is unramified. Furthermore, $\kappa_{0}$ is a  closed immersion, then Corollary \ref{pecigf0ecnr} shows us that $\kappa$ is a pseudo-closed immersion.
\end{proof}

As a consequence of the last result we obtain the following  local description for \'etale morphisms. 

\begin{thm} \label{tppalet}
Let $f\colon\FX \to \FY$ be a morphism in $\sfn$ \'etale at $x \in \FX$. Then there exists an open subset $\FU \subset \FX$ with $x \in \FU$ such that  $f|_{\FU}$ factors as
\[
\FU \xto{\kappa} \FX' \xto{f'} \FY
\]
where $\kappa$ is a completion morphism and $f'$ is an \emph{adic} \'etale morphism.
\end{thm}

\begin{proof}
By the  last theorem  we have that there exists an open formal subscheme $\FU \subset \FX$ with $x \in \FU$ such that  $f|_{\FU}$ factors as 
\[
\FU \xto{\kappa} \FX' \xto{f'} \FY
\]
where $\kappa$ is a pseudo-closed immersion and $f'$ is an adic \'etale morphism. Since $f|_{\FU}$ is \'etale  and $f'$ is an adic \'etale morphism we have that $\kappa$ is \'etale by \cite[Proposition 2.13]{AJP}. Now, applying  Theorem  \ref{caracmorfcompl} it follows that $\kappa$ is a completion morphism.
\end{proof}

\begin{thm} \label{tppall}
Let $f\colon\FX \to \FY$ be a morphism in $\sfn$ smooth at $x \in \FX$. Then there exists an open subset $\FU \subset \FX$ with $x \in \FU$ such that  $f|_{\FU}$ factors as
\[
\FU \xto{\kappa} \FX' \xto{f'} \FY
\]
where $\kappa$ is a completion morphism and $f'$ is an \emph{adic} smooth  morphism.
\end{thm}

\begin{proof}
By Proposition \ref{factpl} there exists an open formal subscheme $\FV \subset \FX$ with $x \in \FV$ such that  $f|_{\FV}$ factors as
\[
\FV \xto{g} \BA^{n}_{\FY} \xto{p} \FY
\]
where $g$ is \'etale  and $p$ is the canonical projection. Applying  the last Theorem to the morphism $g$ we conclude that there exists an open subset $\FU \subset \FX$ with $x \in \FU$ such that  $f|_{\FU}$  factors as 
\[
\FU \xto{\kappa} \FX' \xto{f''} \BA_{\FY}^{n} \xto{p}\FY
\]
where $\kappa$ is a completion morphism, $f''$ is an adic \'etale morphism and $p$ is the canonical projection, from where it follows that $f'= f'' \circ p$ is adic smooth.
\end{proof}

\begin{rem}
Lipman, Nayak and Sastry note in \cite[pag. 132]{LNS} that this Theorem may simplify some developments related to Cousin complexes and duality on formal schemes. See the final part of Remark 10.3.10 of \emph{loc.~cit.\/}
\end{rem}

\begin{ack} 
We have benefited form conversations on these topics and also on terminology with Joe Lipman, Suresh Nayak and Pramath Sastry. The authors thank the Mathematics department of Purdue University for hospitality and support. 

The diagrams were typeset with Paul Taylor's \texttt{diagrams.sty}.
\end{ack}

\end{document}